%% file: main.tex
\documentclass[11pt,onecolumn]{article}

\usepackage{a4wide}
\usepackage{amsmath,amsthm,epsfig,amssymb,amsbsy}
\usepackage{enumerate}
\usepackage{comment}
\usepackage{algorithm}
\usepackage{amsfonts}       %
\usepackage{dsfont}
\usepackage{enumitem}
\usepackage{amssymb}
\usepackage{mathtools}

\usepackage{algpseudocode}

\newcommand\algorithmicinitialize{\textsc{Initialize:}}
\algnewcommand\Initialize{\item[\algorithmicinitialize]}%

\usepackage{nicefrac}

\usepackage{subcaption}
\usepackage{graphicx}

\usepackage{pgfplots}
\pgfplotsset{compat=1.18}

\usepackage{adjustbox}

\usepackage{pgfplotstable}
\usepackage{multirow}

\usepackage[giveninits=true,maxbibnames=9,maxcitenames=2,backend=biber,url=false, doi=false, isbn=false]{biblatex}
\AtEveryBibitem{\clearlist{language}}
\DeclareMathOperator*{\argmin}{arg\,min}

\DeclareMathOperator*\minimize{minimize}

\DeclareMathOperator\stt{subject\ to}
\DeclareMathOperator\st{s.t.}

\DeclareMathOperator{\gph}{gph}

\DeclareMathOperator{\dom}{dom}

\DeclareMathOperator{\relint}{relint}

\DeclareMathOperator{\id}{I}

\usepackage{booktabs,dcolumn,caption}
\newcolumntype{d}[1]{D{.}{.}{#1}} %

\newcommand{\panoc}{\textsc{panoc}}
\newcommand{\Panoc}{\textsc{Panoc}}
\newcommand{\zerofpr}{\textsc{z}ero\textsc{fpr}}

\newcommand{\bR}{\mathbb{R}}

\newcommand{\bN}{\mathbb{N}}

\newcommand{\exR}{\overline{\mathbb{R}}}

\newcommand{\cN}{\mathcal{N}}

\newcommand{\cB}{\mathcal{B}}
\newcommand{\cM}{\mathcal{M}}

\newcommand{\cO}{\mathcal{O}}
\newcommand{\cK}{\mathcal{K}}

\newcommand\downto{\searrow}%

\newcommand\R{\mathbb{R}}
\newcommand\N{\mathbb{N}}
\newcommand\C{\mathcal{C}}
\newcommand\fbe{\varphi_\gamma}
\newcommand\eqdef{:=}
\newcommand\trratio{\rho}

\newcommand\fix{\mathbf{fix}\,}
\newcommand\zer{\mathbf{zer}\,}

\newcommand\seq[1]{(#1)_{k\in\N}}

\makeatletter
\newcommand{\aprox}[3][\@nil]{%
  \def\tmp{#1}%
   \ifx\tmp\@nnil
       \operatorname{prox}_{#3}^{#2}
    \else
         \operatorname{prox}_{#3}^{#1 \star #2}
    \fi}

\newcommand{\aenv}[3][\@nil]{%
  \def\tmp{#1}%
   \ifx\tmp\@nnil
       \operatorname{env}_{#3}^{#2}
    \else
         \operatorname{env}_{#3}^{#1 \star #2}
    \fi}

\newcommand{\bprox}[3][\@nil]{%
  \def\tmp{#1}%
   \ifx\tmp\@nnil
       \operatorname{bprox}_{#3}^{#2}
    \else
        \operatorname{bprox}_{#3}^{#1 #2}
    \fi}
\makeatother

\DeclareMathOperator{\prox}{\mathbf{prox}}
\DeclareMathOperator{\proj}{\Pi}

\usepackage{hyperref}
\hypersetup{
  colorlinks=true,
  linkcolor=blue,
  filecolor=magenta,
  citecolor =magenta,  
  urlcolor=magenta,
  pdftitle={Second-order methods for avoiding nonsmooth saddle points},
}
\usepackage[capitalize]{cleveref}[0.19]

\Crefname{figure}{Figure}{Figures}

\crefformat{equation}{\textup{#2(#1)#3}}
\crefrangeformat{equation}{\textup{#3(#1)#4--#5(#2)#6}}
\crefmultiformat{equation}{\textup{#2(#1)#3}}{ and \textup{#2(#1)#3}}
{, \textup{#2(#1)#3}}{, and \textup{#2(#1)#3}}
\crefrangemultiformat{equation}{\textup{#3(#1)#4--#5(#2)#6}}%
{ and \textup{#3(#1)#4--#5(#2)#6}}{, \textup{#3(#1)#4--#5(#2)#6}}{, and \textup{#3(#1)#4--#5(#2)#6}}

\Crefformat{equation}{#2Equation~\textup{(#1)}#3}
\Crefrangeformat{equation}{Equations~\textup{#3(#1)#4--#5(#2)#6}}
\Crefmultiformat{equation}{Equations~\textup{#2(#1)#3}}{ and \textup{#2(#1)#3}}
{, \textup{#2(#1)#3}}{, and \textup{#2(#1)#3}}
\Crefrangemultiformat{equation}{Equations~\textup{#3(#1)#4--#5(#2)#6}}%
{ and \textup{#3(#1)#4--#5(#2)#6}}{, \textup{#3(#1)#4--#5(#2)#6}}{, and \textup{#3(#1)#4--#5(#2)#6}}

\newtheorem{theorem}{Theorem}[section]
\makeatletter
\newcommand{\settheoremtag}[1]{%
  \let\oldthetheorem\thetheorem%
  \renewcommand{\thetheorem}{#1}%
  \g@addto@macro\endtheorem{%
    \addtocounter{theorem}{-1}%
    \global\let\thetheorem\oldthetheorem}%
  }
\makeatother

\newtheorem{corollary}[theorem]{Corollary}
\newtheorem{lemma}[theorem]{Lemma}
\newlist{lemenum}{enumerate}{1} %
\setlist[lemenum]{label=(\roman*), ref=\thelemma(\roman*), font=\rm}
\crefalias{lemenumi}{lemma}

\crefname{lemma}{Lemma}{Lemmas}
\Crefname{lemma}{Lemma}{Lemmas}

\newtheorem{proposition}[theorem]{Proposition}
\newlist{propenum}{enumerate}{1} %
\setlist[propenum]{label=(\roman*), ref=\theproposition(\roman*), font=\rm}
\crefalias{propenumi}{proposition} 

\newlist{theoremenum}{enumerate}{1} %
\setlist[theoremenum]{label=(\roman*), ref=\thetheorem(\roman*), font=\rm}
\crefalias{theoremenumi}{theorem} 

\newtheorem{definition}[theorem]{Definition}

\newtheorem{assumption}{Assumption}
\crefname{assumption}{Assumption}{Assumptions}
\Crefname{assumption}{Assumption}{Assumptions}

\newtheorem{fact}{Fact}
\crefname{fact}{Fact}{Facts}
\Crefname{fact}{Fact}{Facts}

\newlist{factenum}{enumerate}{1} %
\setlist[factenum]{label=(\roman*), ref=\thefact(\roman*), font=\rm}
\crefalias{factenumi}{fact} 

\newtheorem{condition}[theorem]{Condition}
\crefname{condition}{Condition}{Conditions}
\Crefname{condition}{Condition}{Conditions}

\newtheorem{example}[theorem]{Example}

\newtheorem{remark}[theorem]{Remark}

\newlist{assumenum}{enumerate}{1} %
\setlist[assumenum]{label=(\roman*), ref=\theassumption(\roman*), font=\rm}
\crefalias{assumenumi}{assumption}

\usepackage{xcolor}

\title{
  Second-order methods for provably escaping strict saddle points in composite nonconvex and nonsmooth optimization
}
\author{Alexander Bodard\thanks{KU Leuven,
		Department of Electrical Engineering (ESAT-STADIUS),
		Kasteelpark Arenberg 10, 3001 Leuven, Belgium~
		{\tt%
			\href{mailto:alexander.bodard@kuleuven.be,panos.patrinos@esat.kuleuven.be}{alexander.bodard@kuleuven.be,panos.patrinos@esat.kuleuven.be}%
		}
}
\and Masoud Ahookhosh\thanks{
    Department of Mathematics, University of Antwerp, Antwerp, Belgium
    {\tt%
			\href{mailto:masoud.ahookhosh@uantwerp.be}{masoud.ahookhosh@uantwerp.be}%
		}\\
    A.B.\;and P.P.\;are partially supported by the Research Foundation Flanders (FWO) research projects G081222N, G033822N, G0A0920N and the Research Council KU Leuven C1 project with ID C14/24/103. M.A.\;is partially supported by the Research Foundation Flanders (FWO) research project G081222N and UA BOF DocPRO4 projects with ID 46929 and 48996.
}
\and Panagiotis Patrinos\footnotemark[1]}

\begin{document}
\maketitle

\begin{abstract}
  \input{abstract.tex}

\end{abstract}

\section{Introduction}

\input{introduction/main.tex}

\section{Preliminaries} \label{sec:preliminaries}

\input{preliminaries/main.tex}

\section{Second-order properties of FBE} \label{sec:second-order-properties}

\input{fbe/main.tex}

\section{Nonsmooth trust-region methods on FBE} \label{sec:nonsmooth-tr}

\input{fbtr/main.tex}

\section{Proximal gradient methods with curvilinear linesearch} \label{sec:curvilinear}
\input{curvilinear/main.tex}

\section{Preliminary numerical experiments}\label{sec:numerics}
\input{numerics/main.tex}

\vspace{-5mm}
\section{Conclusion}

\input{conclusion/main.tex}

\begin{appendix}

    \section{Well-definedness and global convergence of Algorithm \ref{alg:fbtr}} \label{app:well-defined-global}

\input{appendix/global-convergence.tex}

\end{appendix}

\printbibliography

\end{document}

%% file: abstract.tex
This study introduces two second-order methods designed to provably avoid saddle points in composite nonconvex optimization problems: (i) a nonsmooth trust-region method and (ii) a curvilinear linesearch method. These developments are grounded in the forward-backward envelope (FBE), for which we analyze the local second-order differentiability around critical points and establish a novel equivalence between its second-order stationary points and those of the original objective. We show that the proposed algorithms converge to second-order stationary points of the FBE under a mild local smoothness condition on the proximal mapping of the nonsmooth term. Notably, for \( \C^2 \)-partly smooth functions, this condition holds under a standard strict complementarity assumption. To the best of our knowledge, these are the first second-order algorithms that provably escape nonsmooth strict saddle points of composite nonconvex optimization, regardless of the initialization. Our preliminary numerical experiments show promising performance of the developed methods, validating our theoretical foundations.

%% file: introduction/main.tex
In this paper, we consider the composite minimization problem
\begin{equation} \label{eq:problem} \tag{P}
    \minimize_{x \in \R^n} \varphi(x) \equiv f(x) + g(x),
\end{equation}
where the function $f$ is smooth but possibly \emph{nonconvex}, and the function $g$ is lower semicontinuous (lsc), proximable, and possibly \emph{nonsmooth}.
The following basic assumptions are used throughout.
\begin{assumption}[Basic assumptions]\label{assump:fbe-basic}
    For the functions defining \eqref{eq:problem}, it holds that
    \begin{assumenum}
        \item \label{assump:fbe-basic-1}$f \in \C^{1, 1}(\R^n)$ (differentiable with $L_f$-Lipschitz continuous gradient);
        \item \label{assump:fbe-basic-2}$g : \R^n \to \overline{\R}$ is proper, lsc, and $\gamma_g$-prox-bounded;
        \item a solution exists, i.e., $\argmin \varphi \neq \emptyset$.
    \end{assumenum}
\end{assumption}\noindent
The composite optimization problem \eqref{eq:problem} covers a wide range of applications, originating in fields such as machine learning, signal processing, statistics and control.
Significant effort has been put into designing efficient solvers that exploit the composite nature of the problem.
In this regard, the \emph{proximal gradient method} (PGM), also known as forward-backward splitting, is a fundamental algorithm which at each iteration involves the computation of a gradient $\nabla f$, and a proximal mapping $\prox_{g}$.
This pair serves as a first-order oracle for \eqref{eq:problem}.
Various acceleration techniques have been proposed for the PGM, including the well-known accelerated proximal gradient method \cite{beck_fast_2009} which attains an improved global convergence rate through Nesterov-type acceleration.
Another line of work revolves around an exact, real-valued penalty function for \eqref{eq:problem}, known as the \emph{forward-backward envelope} (FBE).
Initially proposed in \cite{patrinos_proximal_2013}, the FBE is defined as
\begin{equation} \label{eq:fbe-def}
    \fbe(x) := \inf_{u \in \R^n} \left\{ f(x) + \langle \nabla f(x), u - x \rangle + g(u) + \frac{1}{2\gamma} \Vert u - x \Vert^2 \right\},
\end{equation}
where $\gamma > 0$ corresponds to a proximal stepsize parameter.
Methods like \zerofpr{} \cite{themelis_forward-backward_2018} and \Panoc{} \cite{stella_simple_2017} use the FBE as a merit function in their linesearch procedure to accelerate the PGM with (approximate) second-order directions, and attain local superlinear convergence rates.
Such directions can either be computed based on available first-order information that includes \(\nabla f\) and \(\prox_{g}\), e.g., with L-BFGS, or based on additional second-order information that includes \( \nabla^2 f \) and possibly \( \partial_C (\prox_{g}) \)\footnote{
    The Clarke-generalized Jacobian \( \partial_C(\cdot) \) is used, since \( \prox_{g} \) is not differentiable in general.
}.
This pair constitutes a \textit{second-order oracle} for \eqref{eq:problem}, and is also used by proximal Newton and trust-region methods; see, e.g., \cite{liu_inexact_2024,ouyang_trust_2024} and the references therein.

While existing methods typically focus on a fast local convergence rate, the primary goal of this work is to design second-order methods that provably converge to \emph{second-order stationary points}\footnote{
    A second-order stationary point \(x^\star\) of a nonsmooth function \( \varphi \) is defined in terms of \textit{second-order subderivatives} of \( \varphi \) (cf.\,\cref{sec:optimality-conditions}).
    If \( \varphi \in \C^2 \), this reduces to the familiar condition \( \nabla^2 \varphi(x^\star) \succeq 0 \) with \( \nabla \varphi(x^\star) = 0 \).
} of the nonsmooth problem \eqref{eq:problem}, i.e., they should avoid \emph{nonsmooth strict saddle points} (cf.\,\cref{def:saddle}).
For trust-region methods, both in the smooth setting and on Riemannian manifolds, such guarantees are well established, see e.g., \cite{conn_trust_2000,absil_trust-region_2007}.
An alternative but less common approach is based on \cite{mccormick_modification_1977,more_use_1979} and incorporates a \emph{curvilinear linesearch} procedure in the algorithm, but requires the objective to be \(\C^2\).
Yet, no trust-region or curvilinear linesearch method for \eqref{eq:problem} exists that is known to converge to second-order stationary points for nonsmooth functions \(\varphi\).

To address this open problem, we first investigate novel second-order properties of the FBE under \cref{assump:fbe-basic}. 
Thereafter, we design nonsmooth second-order methods under the following slightly restricted variant of \cref{assump:fbe-basic-2} which guarantees that \(\partial_C(\prox_g)\) is well-defined and nonempty along the iterates.
\begin{assumption}\label{assump:weak-convexity}
    The function \(g : \R^n \to \exR\) is proper, lsc, and \(\rho\)-weakly convex.
\end{assumption}

\subsection{Contributions}
\input{introduction/contributions.tex}

\subsection{Motivation and related work} \label{sec:related-work}
\input{introduction/related-work.tex}

%% file: introduction/contributions.tex
The main contributions of this work can be summarized as follows:
\begin{enumerate}
    \item First, \cref{sec:second-order-properties} studies \textbf{novel second-order properties of the FBE}. 
    \Cref{sec:fbe-local-twice} presents conditions under which $\fbe$ is \emph{of class $\C^2$ around critical points} (\cref{prop:fbe-twice-diff}).
    This nicely complements \cite[Theorem 4.10]{themelis_forward-backward_2018}, which describes the Hessian $\nabla^2 \fbe$, albeit only at critical points. 
    Then, \cref{sec:equivalence} establishes an \emph{equivalence between the second-order stationary points of $\varphi$ and $\fbe$}, a result which appears instrumental for designing proximal gradient-based methods that provably converge to second-order stationary points.
    Contrary to existing works, neither of these properties requires local smoothness of the objective around the given critical point.
    The relevance of both results is further highlighted by the fact that they hold for the broad class of $\C^2$-partly smooth functions $g$ under a strict complementarity condition \eqref{eq:strict-complementarity}.
    \item Second, \cref{sec:nonsmooth-tr} introduces a \textbf{nonsmooth trust-region method} for solving \eqref{eq:problem} which, to the best of our knowledge, is the first to escape \emph{nonsmooth strict saddle points}, irrespective of the initialization.
    The proposed method exploits the equivalence between second-order stationary points of \( \varphi \) and \( \fbe \), and uses a Gauss-Newton approximant that makes the algorithm practically applicable. 
    Global convergence and a local superlinear rate of convergence are established under conventional assumptions.
    \item Third, \cref{sec:curvilinear} presents a \textbf{nonsmooth curvilinear linesearch method} for solving \eqref{eq:problem} that generalizes the globally convergent \zerofpr{} scheme \cite{themelis_forward-backward_2018}.
    As before, a Gauss-Newton approximant is used for practical considerations.
    Similar to \cref{alg:fbtr}, convergence to \emph{second-order stationary points} and a local superlinear rate of convergence are established under suitable assumptions.
    We highlight that, unlike traditional curvilinear linesearch methods, \cref{alg:curvilinear-panoc} does not require a level-boundedness assumption.
\end{enumerate}
Finally, \cref{sec:numerics} illustrates the nonsmooth saddle point avoidance on some toy examples, and demonstrates the practical applicability of \cref{alg:fbtr,alg:curvilinear-panoc} on sparse principal component analysis and phase retrieval problems.

%% file: introduction/related-work.tex
In recent years, problem formulations involving nonconvex functions $f$, and, to some extent involving nonconvex functions $g$, have gained popularity. 
This phenomenon is at least partly driven by the noteworthy success of simple first-order methods in solving nonconvex problems.
In the smooth setting, i.e., $g \equiv 0$, this desirable behavior is often attributed to the \emph{strict saddle property} \cite{lee_first-order_2019,ubl_linear_2023}. 
This property, which is generic among smooth functions, entails that every critical point is either a local minimizer or a strict saddle point, or equivalently, that every second-order stationary point is a local minimizer.
Various machine learning problems are known to possess an even stronger version of the property: in addition to all saddle points being strict, all local minima are also (nearly) global minima \cite{bhojanapalli_global_2016,ge_matrix_2016,ge_no_2017,sun_when_2016,sun_geometric_2018,soltanolkotabi_theoretical_2019,cheridito_gradient_2024, zheng_benign_2024}.
In this context, some authors speak of an \emph{amenable geometry} of the optimization landscape, or of a \emph{benign landscape}.
Interestingly, many first-order methods avoid strict saddles for almost all initializations, and therefore converge to (local) minimizers of problems with the aforementioned properties, for almost any initialization \cite{lee_first-order_2019}.
Second-order methods such as trust-region methods possess similar guarantees, yet irrespective of the initial iterate \cite{conn_trust_2000}.
Recently, Davis et al.\,\cite{davis_escaping_2022,davis_proximal_2022} observed that in a nonsmooth setting, not only the strictly negative curvature, but also the way in which the nonsmoothness manifests itself around the strict saddle point is important.
They showed that simple proximal algorithms converge to local minimizers for almost all initializations under a generalized strict saddle property for nonsmooth functions, which is generic among weakly convex functions.

The convergence of second-order methods to second-order stationary points has been thoroughly studied for the problem of minimizing a smooth objective function.
Trust-region methods arguably constitute the most popular class of methods with such guarantees, see e.g., the monograph \cite{conn_trust_2000}, but also curvilinear linesearch methods possess this property \cite{mccormick_modification_1977,more_use_1979,ferris_nonmonotone_1996}.
Similar results have been obtained for trust-region methods subject to convex constraints \cite[\S 12.4]{conn_trust_2000}, and more recently, trust-region variants have been developed and analyzed for general equality constrained problems, where in particular the smooth objective is minimized over a smooth Riemannian manifold, see e.g., \cite{absil_trust-region_2007,goyens_riemannian_2024}.
Yet, regarding the minimization of a composite objective involving a general nonsmooth term, the existing literature remains limited.
A number of proximal trust-region variants \cite{aravkin_proximal_2022,bodard_pantr_2023,baraldi_proximal_2023,ouyang_trust_2024,zhao_proximal_2024} have been proposed for solving \eqref{eq:problem}, yet none of these comes with second-order stationarity guarantees that mimic those of their smooth counterparts.

%% file: preliminaries/main.tex
\subsection{Notation and basic definitions}
\input{introduction/notation.tex}

\subsection{Forward-backward envelope and proximal gradient method}

The proximal gradient method (PGM) iteratively performs updates
\begin{equation} \label{eq:fbs} \tag{PG}
    x^+ \in T_\gamma(x) := \prox_{\gamma g}(x - \gamma \nabla f(x)),
\end{equation}
for a stepsize $\gamma > 0$.
We define the corresponding \emph{fixed-point residual} as $R_\gamma(x) := \gamma^{-1} (x - T_\gamma(x))$
and say that a point $x^\star \in \fix T_\gamma$ is \emph{critical} if $x^\star \in T_\gamma(x^\star)$ for some $\gamma \in (0, \gamma_g)$.
The limit points of \eqref{eq:fbs} are critical points when $\gamma \in (0, \min \{\nicefrac{1}{L_f}, \gamma_g\})$, and this independently of the choice of $x^+ \in T_\gamma(x)$ \cite[Lemma 2]{bolte_proximal_2014}.
The FBE \eqref{eq:fbe-def}, which has been used to accelerate PGM with fast Newtonian directions, can also be expressed more explicitly in terms of the Moreau envelope as
\begin{equation} \label{eq:fbe-moreau-formulation}
    \fbe(x) = f(x) - \frac{\gamma}{2} \Vert \nabla f(x) \Vert^2 + g^\gamma(x - \gamma \nabla f(x)).
\end{equation}
Thus, it is clear that $\fbe$ inherits the regularity properties of the Moreau envelope $g^\gamma$, including that for any $\gamma \in (0, \gamma_g)$, the function $\fbe$ is real-valued and strictly continuous on $\R^n$ \cite[Proposition 4.2]{themelis_forward-backward_2018}.
Other useful properties are given below.

\begin{fact} [{\cite[Prop.\,4.3 and Th.\,4.4]{themelis_forward-backward_2018}}] \label{prop:fbe-decrease-minimizers}
    Let \( \gamma \in (0, \gamma_g) \) be fixed. Then, the following statements hold:
    \begin{factenum}
        \item \( \fbe \leq \varphi \);
        \item \( \varphi(\bar x) \leq \fbe(x) - \frac{1 - \gamma L_f}{2\gamma} \Vert x - \bar x \Vert^2 \) for all \( x \in \R^n \) and \( \bar x \in T_\gamma(x) \);
        \item $\varphi(x^\star)=\fbe(x^\star)$ for all $x^\star \in \fix T_\gamma$;
        \item if $\gamma \in \left(0, \min \left\{ 1 / L_f, \gamma_g \right\} \right)$, then $\inf \varphi=\inf \fbe$ and $\argmin \varphi = \argmin \fbe$.
    \end{factenum}
\end{fact}

The local differentiability of $\fbe$ around a critical point $x^\star$ is closely related to the local differentiability of $g^\gamma$ around the forward point $x^\star - \gamma \nabla f(x^\star)$, as shown by \eqref{eq:fbe-moreau-formulation}, and this is in turn related to \emph{prox-regularity} of $g$ at $x^\star$ for $-\nabla f(x^\star)$, as defined below.
\begin{definition}[Prox-regularity \cite{poliquin_prox-regular_1996}]
    A proper, lsc function $g : \R^n \to \exR$ is \emph{prox-regular} at $\bar{x}$ for $\bar{v} \in \partial g(\bar{x})$ if there exist $r \geq 0, \epsilon > 0$ such that for all $x' \in \cB(\bar{x}; \epsilon)$ and
    \begin{equation} \label{eq:prox-regularity-xv}
        (x, v) \in \gph \partial g \quad \st \quad x \in \cB(\bar{x}; \epsilon), \quad v \in \cB(\bar{v}; \epsilon), \quad g(x) \leq g(\bar{x}) + \epsilon,
    \end{equation}
    it holds that $g(x') \geq g(x) + \langle v, x' - x \rangle - \frac{r}{2} \Vert x' - x \Vert^2$.
\end{definition} \noindent
The prox-regularity of $g$ is a mild requirement enjoyed by a wide class of functions that includes all strongly amenable functions \cite[Proposition 13.32]{rockafellar_variational_1998}. 
Taking $\epsilon = +\infty$ and $r = \rho$, it is clear that $\rho$-weakly convex functions satisfy this requirement globally for any subgradient.
Thus, prox-regularity can be thought of as a ``local" version of weak-convexity, in the sense that this subgradient inequality holds for all $(x, v) \in \gph \partial g$ that are sufficiently close to $(\bar{x}, \bar{v})$, with the additional restriction that $g(x)$ is not too different from $g(\bar{x})$.
This concept is often formalized by the so-called $g$-attentive $\epsilon$-localization of $\partial g$ at $(\bar{x}, \bar{v})$, defined as the mapping
\begin{equation*}
    \widetilde T_g (x) = \begin{cases}
        \left\{ v \in \partial g(x) \mid v \in \cB(\bar{v}; \epsilon) \right\} & \text{ if } x \in \cB(\bar{x}; \epsilon) \text{ and } g(x) \leq g(\bar{x}) + \epsilon,\\
        \emptyset & \text{ otherwise.}
    \end{cases}
\end{equation*}
Observe that we can thus replace \eqref{eq:prox-regularity-xv} in the definition of prox-regularity by $(x, v) \in \gph \widetilde T_g$.
Moreover, we note that when $g$ is assumed to be subdifferentially continuous, as in \cite{hang_smoothness_2023}, the $g$-attentiveness condition $g(x) \leq g(\bar{x}) + \epsilon$ becomes redundant.

\begin{fact}[{\cite[Theorem 4.7]{themelis_forward-backward_2018}}] \label{prop:fbe-locally-c1}
    Suppose that $f$ is of class $\C^2$ around a critical point $x^\star$, and that $g$ is prox-regular at $x^\star$ for $-\nabla f(x^\star)$.
    Then, for all $\gamma > 0$ sufficiently small there exists a neighborhood $U_{x^\star}$ of $x^\star$ on which the following hold:
    \begin{factenum}
        \item \( T_\gamma \) and \( R_\gamma \) are strictly continuous, and in particular single-valued;
        \item \( \fbe \in \C^1 \) with \( \nabla \fbe(x)=Q_\gamma(x) R_\gamma(x) \), where \(Q_\gamma(x)=I-\gamma \nabla^2 f(x)\).
    \end{factenum}
\end{fact}

\subsection{Optimality conditions and (strict) saddle points} \label{sec:optimality-conditions}

A generalization of Fermat's rule for nonsmooth functions states that a local minimizer \(x^\star \in \dom \varphi\) of a proper function \(\varphi : \R^n \to \exR\) satisfies \(0 \in \partial \varphi(x^\star)\) \cite[Thm.\,10.1]{rockafellar_variational_1998}. 
Moreover, the second-order necessary conditions of optimality can be extended to nonsmooth functions.
\begin{fact}[{\cite[Theorem 13.24]{rockafellar_variational_1998}}] \label{prop:second-order-conditions}
    If \(x^\star \in \dom \varphi\) is a local minimizer of a proper function \(\varphi : \R^n \to \exR\), then \(0 \in \partial \varphi(x^\star)\) and $\mathrm d^2 \varphi(x^\star, 0)[d] \geq 0$ for all $d \in \bR^n$.
\end{fact}\noindent
When minimizing general nonconvex problems, only convergence of the iterates to \emph{stationary points} can be guaranteed, i.e., to points satisfying a necessary condition of optimality.
\begin{definition}[Stationary points] \label{def:stationarity}
    Let \(\varphi : \R^n \to \exR\). We say that a point \(x^\star \in \dom \varphi\) is
    \begin{enumerate}[label=(\roman*), font=\rm]
        \item a first-order stationary point of \(\varphi\) if \(x^\star \in \zer \partial \varphi\);
        \item a second-order stationary point of \(\varphi\) if \(x^\star \in \zer \partial \varphi\) and $\mathrm d^2 \varphi(x^\star, 0)[d] \geq 0$ for all $d \in \bR^n$.
    \end{enumerate}
\end{definition}\noindent
If \( \varphi \in \C^2 \), then the definition of a second-order stationary point reduces to the classical condition \( \nabla^2 \varphi(x^\star) \succeq 0 \) with \( \nabla \varphi(x^\star) = 0 \).
Note that the limit points of PGM iterations are critical points, which in turn implies first-order stationarity.
In turn, criticality is a necessary condition of optimality for problems of the form \eqref{eq:problem}.

\begin{fact}[{\cite[Proposition 3.5]{themelis_forward-backward_2018}}]
    Under \cref{assump:fbe-basic} the following statements hold:
    \begin{factenum}
        \item (criticality $\Rightarrow$ first-order stationarity) $\fix T_\gamma \subseteq \zer \partial \varphi$ for all $\gamma \in (0, \gamma_g)$;
        \item (optimality $\Rightarrow$ criticality) $\argmin \varphi \subseteq \fix T_\gamma$ for all $\gamma$ sufficiently small.
    \end{factenum}
\end{fact}\noindent
Finally, we define (strict) saddle points of a nonsmooth function in terms of the first-order and second-order stationarity measures.
\begin{definition}[(Strict) saddle points] \label{def:saddle}
    Let \(\varphi : \R^n \to \exR\). We say that a first-order stationary point \(x^\star \in \dom \varphi\) is
    \begin{enumerate}[label=(\roman*), font=\rm]
        \item a \emph{saddle point} if \(x^\star\) is neither a local minimizer, nor a local maximizer;
        \item a \emph{strict saddle point} if \(x^\star\) is not a second-order stationary point.
    \end{enumerate}
\end{definition}\noindent
By \cref{def:stationarity}, this means that a strict saddle point \(x^\star \in \zer \partial \varphi\) satisfies $\mathrm d^2 \varphi(x^\star, 0)[d] < 0$ for some $d \in \bR^n$.
If \(\varphi \in \C^2\), then we recover the classical condition \(\nabla \varphi(x^\star) = 0\) and \(\langle \nabla^2 \varphi(x^\star) d, d \rangle < 0\) for some \(d \in \bR^n\).
If \(\varphi\) is not of class \(\C^2\) locally around a saddle point \(x^\star\), then we say that \(x^\star\) is a \emph{nonsmooth saddle point}.
We remark that the work \cite{liu_envelope_2019} also considers problems of the form \eqref{eq:problem}, but assumes \(\varphi\) to be smooth locally around critical points.
In contrast, this work focuses on nonsmooth saddle points, which arise naturally in the context of \eqref{eq:problem}.
The following examples illustrate this.
\begin{example}[Nonconvex quadratic subject to box constraints] \label{ex:quadratics-box}
    \begin{figure}
        \begin{subfigure}{0.5\textwidth}
            \centering
            \includegraphics[width=\textwidth]{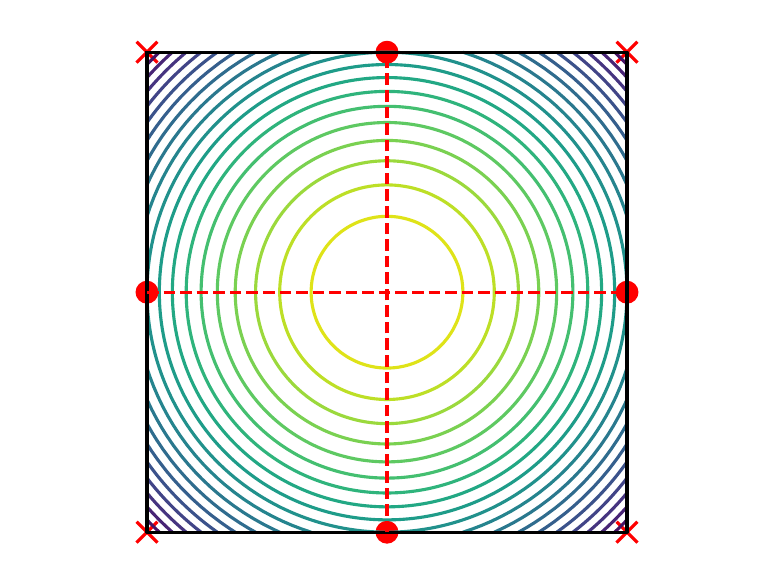}
            \captionsetup{width=1\textwidth}
            \caption[]{Level curves of $\varphi(x, y) = -x^2 - y^2 + \delta_{[-1, 1]}(x, y)$.}
            \label{fig:nonconvex-quadratic}
        \end{subfigure}
        \begin{subfigure}{0.49\textwidth}
            \centering
            \includegraphics[width=\textwidth]{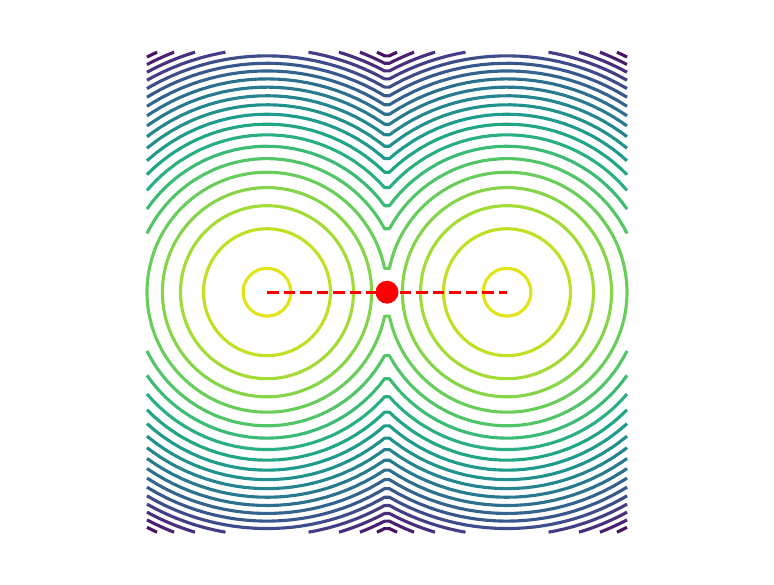}
            \captionsetup{width=.85\textwidth}
            \caption[]{Level curves of $\varphi(x, y) = - x^2 - y^2 + \vert x \vert$.}
            \label{fig:nonconvex-sparsity}
        \end{subfigure}
        \caption{Crosses indicate minimizers, dots correspond to saddle points. 
        The objective $\varphi$ is \emph{nonsmooth} at all minimizers and saddle points.
        The PGM escapes the saddle points, unless it is initialized on the dashed lines.
        }
        \label{fig:motivating-examples}
    \end{figure}
    Consider a problem where $f$ is a nonconvex quadratic function and $g = \delta_C$ is the indicator of a set of box constraints $C$.
    In a (local) minimizer where at least some of the box constraints are active, the function $g = \delta_C$ is discontinuous.
    \Cref{fig:nonconvex-quadratic} depicts this situation, and visualizes how \emph{nonsmooth saddle points} naturally arise at boundary points.
\end{example}
\begin{example}[Sparsity by $\ell_1$-regularization] \label{ex:l1}
    Sparsity is often enforced through a regularization term $g(x) = \Vert x \Vert_1$.
    In this setting, one expects -- or at least hopes -- to find a critical point with a number of zero components.
    However, around such points $g$ is not $\C^1$.
    \Cref{fig:nonconvex-sparsity} visualizes a \emph{nonsmooth saddle point} that arises in this way when $f$ is a nonconvex function.
\end{example}

\subsection{Partial smoothness} \label{sec:partly-smooth}
\input{fbe/partial-smoothness.tex}

%% file: introduction/notation.tex
We denote by $\exR = \bR \cup \left\{ +\infty \right\}$ the extended real line, and by \( \id_n \) the $n \times n$ identity operator.
If the dimension \( n \) is clear from the context, we simply write $\id$. 
Unless mentioned otherwise, \( \langle \cdot, \cdot \rangle \) denotes the standard Euclidean inner product, and \( \Vert \cdot \Vert \) the induced norm.
The smallest eigenvalue of a symmetric matrix \( B \in \R^{n \times n} \) is referenced by \( \lambda_{\min}(B) \), and \( \relint C \) denotes the relative interior of a set \( C \subseteq \R^n \).
The open and closed balls of radius $r \geq 0$ centered in $x \in \R^n$ are denoted by $\cB(x;r)$ and \(\overline\cB(x;r)\), respectively.
By a slight abuse of terminology, we call a sequence \( \seq{x^k} \) summable if \( \sum_{k \in \N} \Vert x^k \Vert \) is finite, and square-summable if \( \seq{\Vert x^k \Vert^2} \) is summable.

We say that a function \(f : \R^n \to \exR\) is \(\rho\)-weakly convex if \(x \mapsto f + \frac{\rho}{2} \Vert x \Vert\) is convex.
The set of limit points of \( \seq{x^k} \) is denoted by \( \omega(\seq{x^k}) \).
We follow the definitions of \cite{rockafellar_variational_1998} in terms of strict continuity and (strict) differentiability of functions and mappings.
The set of \( k \) times continuously differentiable functions $f : \R^n \to \R$ is denoted by \( \C^k \), the set of functions \( f \in \C^k \) with locally Lipschitz continuous \( k \)'th derivative is denoted by \( \C^{k+} \), and the set of functions \( f \in \C^1 \) with Lipschitz continuous gradient is denoted by \( \C^{1, 1} \).
We denote by \( \mathrm d g(\bar x)(\bar w) \) the (lower) subderivative of \( g : \R^n \to \exR \) at \( \bar x \in \dom g \) for \( \bar w \) \cite[Def.\,8.1]{rockafellar_variational_1998}, and by \( \mathrm d^2 g(\bar x, v)[d] \) the second subderivative at \( \bar x \in \dom g \) for \( v \in \R^n \) and \( w \in \R^n \) \cite[Def.\,13.3]{rockafellar_variational_1998} where \( \dom g := \{ x\in \R^n \mid g(x) < \infty \} \) is the effective domain of \( g \).
A vector \( v \in \partial g(x) \) is a subgradient of a proper function \( g : \R^n \to \exR \) in \( x \in \R^n \), where the subdifferential \( \partial g(x) \) is defined as in \cite[Def.\,8.3]{rockafellar_variational_1998}, that is,
\begin{equation*}
    \partial g(x) = \left\{ v \in \bR^n \mid \exists (x^k)_{k \in \bN} \to x, (v^k \in \hat{\partial}g(x^k))_{k \in \bN} \to v\; \st.\; g(x^k) \to g(x) \right\},
\end{equation*}
where $\hat\partial g(x)$ denotes the set of \emph{regular subgradients} of $g$ at $x$ as defined in \cite[Def.\,8.3]{rockafellar_variational_1998}.

For a given stepsize parameter $\gamma > 0$, the \emph{Moreau envelope} $g^\gamma$ and the \emph{proximal mapping} or proximal operator $\prox_{\gamma g}$ of $g$ are defined as
\begin{equation*}
    \begin{aligned}
        g^\gamma(x) \eqdef \inf_{z} \left\{ g(z) + \frac{1}{2\gamma} \Vert z - x \Vert^2 \right\}, \quad
        \prox_{\gamma g}(x) \eqdef \argmin_z \left\{ g(z) + \frac{1}{2\gamma} \Vert z - x \Vert^2 \right\}.
    \end{aligned}
\end{equation*}
The function \( g \) is \( \gamma_g \)-prox-bounded if for any \( \gamma \in (0, \gamma_g) \) the function \( g + \frac{1}{2\gamma} \Vert \cdot \Vert^2 \) is bounded below on \( \R^n \).
For any such \( \gamma \), the proximal mapping \( \prox_{\gamma g} \) is nonempty and compact-valued, and the Moreau envelope is finite \cite[Thm.\,1.25]{rockafellar_variational_1998}.
Given a nonempty and closed set \( S \subseteq \R^n \), we define its indicator function \( \delta_S : \R^n \to \exR \) as \( \delta_S(x) = 0 \) if \(x \in S\) and \( \delta_S(x) = \infty \) otherwise, and denote by \( \proj_S : \R^n \rightrightarrows \R^n \) the (set-valued) projection onto \( S \).
For a set-valued mapping \( T : \R^n \rightrightarrows \R^n \), the graph of \( T \) is \( \gph T = \{ (x, y) \mid y \in T(x) \} \), the zeros of \( T \) are \( \zer T = \{ x \mid 0 \in T(x) \} \) and the fixed-points of \( T \) are \( \fix T = \{ x \mid x \in T(x) \} \).

%% file: fbe/partial-smoothness.tex
The class of $\C^2$-partly smooth functions was initially introduced in \cite{lewis_active_2002} and further analyzed in \cite{mifflin_primal-dual_2003,hare_nonsmooth_2004,hare_identifying_2004,drusvyatskiy_optimality_2014,daniilidis_orthogonal_2014,lewis_partial_2013,hang_fresh_2024}.
Locally, these nonsmooth functions are associated with a smooth manifold, often called \emph{active manifold}, and thus possess a certain smooth structure.
Recall that a set $\cM \subseteq \R^n$ is a $\C^2$-smooth manifold of codimension $m \leq n$ around $\bar x \in \R^n$ if $\bar x \in \cM$ and if there exists an open neighborhood $\cO \subseteq \R^n$ containing $\bar x$ and a $\C^2$ mapping $\Phi : \cO \to \R^m$ for which
\begin{equation*}
    \cM \cap \cO = \left\{ x \in \cO \mid \Phi(x) = 0 \right\} \quad \text{and} \quad \nabla \Phi(\bar x) \text{ is surjective.}
\end{equation*}
However, these functions behave ``sharply'' along normal directions of the active manifold.
\begin{definition}[Partial smoothness {\cite{lewis_active_2002}}] \label{def:partial-smoothness}
    Suppose that the set $\cM \subseteq \bR^n$ contains the point $\bar x$. The function $g$ is $\C^2$-\emph{partly} smooth at $\bar x$ relative to $\cM$ if $\cM$ is a $\C^2$-smooth manifold around $\bar x$ and the following properties hold:
    \begin{enumerate}[label=(\roman*), font=\rm]
        \item[(i)] (restricted smoothness) the restriction $g_{\mid \cM} := g + \delta_{\cM}$ is $\C^2$ around $\bar x$;
        \item[(ii)] (regularity) at every point $x \in \cM$ close to $\bar x$, the function $g$ is subdifferentially regular in the sense of \cite[Def.\,7.25]{rockafellar_variational_1998} and has nonempty subdifferential $\partial g(x)$;
        \item[(iii)] (normal sharpness) $\mathrm d g(\bar x)(-w) > - \mathrm d g(\bar x) (w)$ for all normal directions $0 \neq w \in N_{\cM}(\bar x)$;
        \item[(iv)] (subdifferential continuity) the subdifferential map $\partial g$ is continuous at $\bar x$ relative to $\cM$.
    \end{enumerate} 
\end{definition}\noindent
\begin{example}
\Cref{fig:partial-smoothness} illustrates \cref{def:partial-smoothness} for the function \(g(x, y) = \vert x \vert + \frac{1}{2} y^2\), which is \(\C^2\)-partly smooth at \(\bar x =(0, 0)\) relative to \(\cM := \{(0, y) \mid y \in \R\}\).
Indeed, the restricted smoothness condition is easily verified, and the normal sharpness follows from the fact that for any \(w = (\bar w_x, 0) \in N_{\cM}(0, 0) = \{(w_x, 0) \mid w_x \in \R\}\) we have \(\mathrm d g(0, 0)(w) = \vert \bar w_x \vert\).
Finally, the regularity and subdifferential continuity conditions follow from the convexity of \(g\).
\begin{figure}
    \centering
    \includegraphics[width=\textwidth]{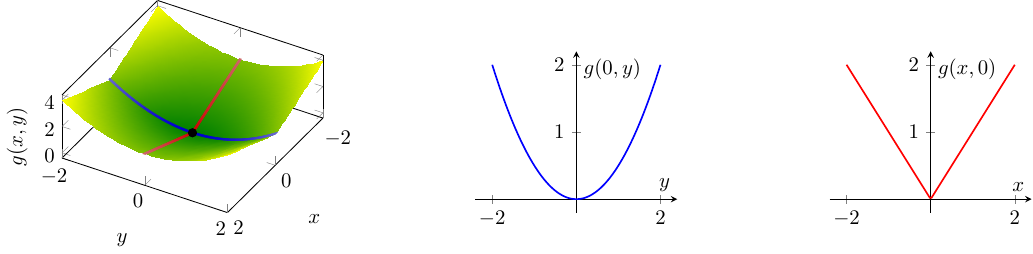}
    \captionsetup{width=\textwidth}
    \caption{
        The function \(g(x, y) = \vert x \vert + \frac{1}{2} y^2\) (left) is \(\C^2\)-partly smooth at the origin relative to \(\cM := \{(0, y) \mid y \in \R \}\).
        Around \((0, 0)\), the restriction of \(g\) to \(\cM\) is \(\C^2\) smooth (middle), and \(g\) behaves ``sharply'' along normal directions \(w \in N_{\cM}(0, 0) = \{ (w_x, 0) \mid w_x \in \R \}\) (right).}
    \label{fig:partial-smoothness}
\end{figure}
\end{example}
Many functions \( g \) that are commonly used in the context of \eqref{eq:problem} are \( \C^2 \)-partly smooth on their domain.
In particular, we highlight several popular regularizers, including the sparsity-promoting $\ell_1$-norm, the $\ell_{1,2}$-norm, the classical $\ell_2$-norm, the anti-sparsity promoting $\ell_\infty$-norm, the low rank-promoting nuclear norm, the total variation semi-norm, the $\ell_0$ pseudo-norm, and the rank function \cite[\S 5.2]{liang_convergence_2016}.
In addition, functions of which the epigraph is a polyhedral set, are $\C^2$-partly smooth at any point in their domain \cite[Ex.\,3.4]{lewis_active_2002}.

Various other smooth substructures of nonsmooth functions have been related to the notion of partial smoothness \cite{hare_functions_2006}.
    For example, a function $h : \R^n \to \exR$ is $\C^2$-decomposable \cite{shapiro_class_2003} at a point $\bar x \in \R^n$ if $h(\bar x)$ is finite and $h$ can be locally represented as
    \begin{equation*}
        h(x) = h(\bar x) + \theta(F(x)), \quad \text{for} \quad x \in \cO,
    \end{equation*}
    where $\cO \subseteq \R^n$ is an open neighborhood around $\bar x$, $\theta: \R^m : \to \exR$ is a proper, lsc, and sublinear function, and $F : \R^n \to \R^m$ is a $\C^2$-smooth mapping with $F(\bar x) = 0$.
    $\C^2$-decomposable functions are $\C^2$-partly smooth under a nondegeneracy or transversality assumption \cite{shapiro_class_2003,hare_functions_2006}. 
    These nondegenerate $\C^2$-decomposable functions are also referred to as \emph{reliably $\C^2$-decomposable} \cite{hang_smoothness_2023}.
    Examples include polyhedral functions and indicators of $\C^2$-cone reducible sets\footnote{We refer to \cite[Def.\,3.135]{bonnans_perturbation_2000} for a formal definition and additional examples.}\cite[Ex.\,5.3a]{hang_smoothness_2023}.
    At all of their points, polyhedral convex sets \cite[Ex.\,3.139]{bonnans_perturbation_2000}, the second-order cone, and the positive semidefinite cone \cite[Ex.\,3.140]{bonnans_perturbation_2000} are $\C^2$-cone reducible, and hence indicators of these sets are $\C^2$-partly smooth.

%% file: fbe/main.tex
This section presents novel second-order properties of the forward-backward envelope that are fundamental in establishing second-order convergence of the methods presented in the subsequent sections.
First, \cref{sec:fbe-local-twice} shows that $\fbe$ is of class \( \C^2 \) locally around a critical point $x^\star$ if the proximal mapping \( \prox_{\gamma g} \) with stepsize \( \gamma > 0 \) is of class \( \C^1 \) locally around the forward point \( x^\star - \gamma \nabla f(x^\star) \).
This extends existing results that describe the differentiability of $\fbe$ at a critical point $x^\star$.
We also derive an equivalent characterization in terms of the second-order epi-derivatives of \( g \), and show that this condition holds for $\C^2$-partly smooth functions $g$ under a mild strict complementarity condition.
Finally, \cref{sec:equivalence} derives an equivalence between the second-order stationary points of $\varphi$ and $\fbe$ under mild conditions.

\subsection{Twice differentiability around critical points} \label{sec:fbe-local-twice}
\input{fbe/local-C2.tex}

\subsection{Equivalence of second-order stationary points} \label{sec:equivalence}
\input{fbe/equivalence-strict-saddles.tex}

%% file: fbe/local-C2.tex
By expressing $\fbe$ as in \eqref{eq:fbe-moreau-formulation}, we observe that the properties of $\fbe$ and $g^\gamma$ are closely related.
In particular, local twice differentiability of $\fbe$ around a critical point $x^\star$ requires local twice differentiability of $g^\gamma$ around the forward point $x^\star - \gamma \nabla f(x^\star)$, which is in turn equivalent to \( \prox_{\gamma g} \) being of class \( \C^1 \) around that forward point.
Poliquin and Rockafellar \cite{poliquin_generalized_1996} describe these and other equivalent characterizations, one of which involves prox-regularity of $g$ and twice epi-differentiability of $g$ \cite[Def.\,13.6]{rockafellar_variational_1998}.
Remark that strict twice epi-differentiability of $g$ at $x^\star$ implies that the second-order epi-derivative of $g$ is generalized quadratic, as can be seen from \cite[Th 4.1]{poliquin_generalized_1996}. 
\begin{assumption} \label{assump:fbe-twice-epi-all}
    With respect to a given critical point $x^\star \in \R^n$, the function $g$ is prox-regular at $x^\star$ for $-\nabla f(x^\star)$, and
    \begin{assumenum}
        \item \label{assump:fbe-twice-epi} twice epi-differentiable at $x^\star$ for $-\nabla f(x^\star)$, with its second-order epi-derivative being generalized quadratic;
        \item \label{assump:fbe-twice-epi-strict} \emph{strictly} twice epi-differentiable at $x^\star$ for $-\nabla f(x^\star)$;
        \item \label{assumption:fbe-local-twice-diff} strictly twice epi-differentiable at $x$ for $v$ with respect to all $(x, v) \in \gph \widetilde T_g$ sufficiently near $(x^\star, -\nabla f(x^\star))$, where $\widetilde T_g$ is the $g$-attentive $\epsilon$-localization of $\partial g$ at $(x^\star, -\nabla f(x^\star))$.
    \end{assumenum}
\end{assumption}\noindent

The theorems in \cite{poliquin_generalized_1996} are formulated under simplifying assumptions, e.g., twice differentiability of $g^\gamma$ around \textit{global minimizers} of $g$ is studied. 
This may seem restrictive for our purposes since we are interested in the second-order differentiability of $g^\gamma$ around $x^\star - \gamma \nabla f(x^\star)$ where $x^\star$ is a critical point.
Yet, by appropriately tilting $g$ and by adding a sufficiently large quadratic term, the results from Poliquin and Rockafellar can be translated to our setting.  
This trick of considering a `tilted' function was implicitly used by \cite{stella_forwardbackward_2017,themelis_forward-backward_2018} to describe twice differentiability of the FBE at critical points (cf.\,\cite[Theorem 4.10]{themelis_forward-backward_2018}).
For the sake of completeness, we describe and prove a variant of \cite[Theorem 3.8, 4.1, 4.4]{poliquin_generalized_1996}, which relates \cref{assump:fbe-twice-epi-all} to differentiability of $\prox_{\gamma g}$ in our setting.

\begin{lemma} \label{th:prox-diff}
    Suppose that \cref{assump:fbe-basic} holds.
    For all $\gamma > 0$ sufficiently small, the following statements regarding a critical point $x^\star \in \R^n$ hold:
    \begin{lemenum}
        \item \label{th:prox-diff-point} $\prox_{\gamma g}$ is differentiable at $x^\star - \gamma \nabla f(x^\star)$ iff \cref{assump:fbe-twice-epi} holds with respect to $x^\star$;
        \item \label{th:prox-diff-point-strict} $\prox_{\gamma g}$ is \emph{strictly} differentiable at $x^\star - \gamma \nabla f(x^\star)$ iff \cref{assump:fbe-twice-epi-strict} holds with respect to $x^\star$;
        \item \label{th:prox-diff-neighborhood} $\prox_{\gamma g}$ is continuously differentiable on a neighborhood around the point $x^\star - \gamma \nabla f(x^\star)$ iff \cref{assumption:fbe-local-twice-diff} holds with respect to $x^\star$.
    \end{lemenum}
\end{lemma}
\begin{proof}
    By criticality of $x^\star$, it follows from \cite[Theorem 3.4 (i)]{themelis_forward-backward_2018} that
    \begin{equation*}
        g(x) \geq g(x^\star) + \langle -\nabla f(x^\star), x - x^\star \rangle - \frac{1}{2\mu} \Vert x - x^\star \Vert^2, \quad \forall x \in \R^n,
    \end{equation*}
    for all $\mu > 0$ sufficiently small \cite[Theorem 3.4 (ii)]{themelis_forward-backward_2018}.
    By possibly reducing the stepsize \( \gamma \), we can ensure that the above equation holds for \( \mu = 2 \gamma \).
    This implies that the function
    \begin{equation*}
        \tilde g(x) \eqdef g(x) + \langle \nabla f(x^\star), x - x^\star \rangle + \frac{1}{2\mu} \Vert x - x^\star \Vert^2.
    \end{equation*}
    satisfies $\tilde g (x) \geq \tilde g(x^\star)$ for all \(x \in \R^n\), and therefore we have that $x^\star \in \argmin \tilde g$.
    The function $\tilde g$ is prox-regular at $x^\star$ for $0$ iff $g$ is prox-regular at $x^\star$ for $-\nabla f(x^\star)$ \cite[Exercise 13.35]{rockafellar_variational_1998}.
    Moreover, $\tilde g$ is (strictly) twice epi-differentiable at $x^\star$ for $0$ iff $g$ is (strictly) twice epi-differentiable at $x^\star$ for $-\nabla f(x^\star)$ \cite[Proposition 2.10]{rockafellar_first-_1988}.
    By \cite[Theorem 3.8]{poliquin_generalized_1996} we therefore have that \cref{assump:fbe-twice-epi} holds iff $\prox_{\mu \tilde g}$ is differentiable at $x^\star$ for $\mu > 0$ sufficiently small, and by \cite[Theorem 4.1]{poliquin_generalized_1996} that \cref{assump:fbe-twice-epi-strict} holds iff $\prox_{\mu \tilde g}$ is strictly differentiable at $x^\star$ for $\mu > 0$ small enough.
    The first two claims follow from tilting $\tilde g$ back to $g$.
    In particular, we can substitute
    \begin{equation}\label{prf:prox-substitution}
        \begin{aligned}
            \prox_{\gamma g}(x - \gamma \nabla f(x^\star)) &= \argmin_{z} \left\{ g(z) + \frac{1}{2\gamma} \Vert z - (x - \gamma \nabla f(x^\star)) \Vert^2 \right\}\\
            &= \argmin_z \left\{ g(z) + \langle \nabla f(x^\star), z - x \rangle + \frac{1}{2\gamma} \Vert z - x \Vert^2 \right\}\\
            &= \argmin_z \left\{ \tilde g(z) + \left( \frac{1}{2\gamma} - \frac{1}{2\mu} \right) \Vert z - x \Vert^2 \right\} = \prox_{\mu \tilde g}(x)
        \end{aligned}
    \end{equation}
    in the definition of differentiability of $\prox_{\mu \tilde g}$ at $x^\star$, i.e.,
    \begin{equation*}
        \frac{\prox_{\mu \tilde{g}}(x) - \prox_{\mu \tilde{g}}(x^\star) - P_{\gamma}(x^\star)(x - x^\star)}{\Vert x - x^\star \Vert} \to 0, \quad x \to x^\star, \quad x \neq x^\star,
    \end{equation*}
    which immediately yields differentiability of $\prox_{\gamma g}$ at $x^\star - \gamma \nabla f(x^\star)$, with Jacobian
    \begin{equation*}
        J \prox_{\gamma g}(x^\star - \gamma \nabla f(x^\star)) = P_\gamma(x^\star) \eqdef J \prox_{\mu \tilde g}(x^\star).
    \end{equation*}
    This proves \cref{th:prox-diff-point}.
    In a similar way, \cref{th:prox-diff-point-strict} follows by substitution of \eqref{prf:prox-substitution} in the strict differentiability of $\prox_{\mu \tilde g}$.
    As for the third claim, remark that $\tilde g$ is strictly twice epi-differentiable at $x$ for $v$ w.r.t.\,all $(x, v) \in \gph T_{\tilde g}$ sufficiently near to $(x^\star, 0)$ if \cref{assump:fbe-twice-epi-all} holds \cite[Proposition 2.10]{rockafellar_first-_1988}.
    By \cite[Theorem 4.4]{poliquin_generalized_1996}, this is in turn equivalent to $\prox_{\mu \tilde g}$ being continuously differentiable on a neighborhood $U_{x^\star}$ around $x^\star$ for $\mu = 2\gamma $ sufficiently small.
    By again substituting \eqref{prf:prox-substitution}, it follows that $\prox_{\gamma g}$ is differentiable at $x - \gamma \nabla f(x^\star)$ for $x \in U_{x^\star}$.
\end{proof}

\Cref{th:prox-diff} details necessary and sufficient conditions for $\prox_{\gamma g}$ to be differentiable at and around critical points.
Under appropriate smoothness assumptions on $f$, we can use \cref{th:prox-diff-neighborhood} to describe the Hessian of $\fbe$ \textit{around} critical points, thus complementing the results from \cite[Theorem 4.10]{themelis_forward-backward_2018} which describe the Hessian of $\fbe$ only \textit{at} critical points.
\begin{assumption}\label{assump:fbe-twice-diff}
    With respect to a critical point \( x^\star \in \R^n \) and for a given stepsize \( \gamma > 0 \):
    \begin{assumenum}
        \item \( f \) is of class \( \C^3 \), locally around \( x^\star \);
        \item \( \prox_{\gamma g} \) is of class \( \C^1 \) locally around \( x^\star - \gamma \nabla f(x^\star) \).
    \end{assumenum}
\end{assumption}
\begin{theorem}[Local twice differentiability of $\fbe$] \label{prop:fbe-twice-diff}
    Suppose that \cref{assump:fbe-basic} holds.
    If \cref{assump:fbe-twice-diff} holds with respect to a critical point \( x^\star \in \R^n \) and a sufficiently small \(
    \gamma > 0\),
    then there exists a neighborhood $U_{x^\star}$ of $x^\star$ on which $R_\gamma$ is continuously differentiable and $\fbe \in \C^2$. 
    In particular, for all $x \in U_{x^\star}$ we have
    \begin{equation*}
        \nabla^2 \fbe(x) = Q_\gamma(x) J R_\gamma(x) + M_\gamma(x), \quad \text{ with } \quad J R_\gamma(x) = \frac{1}{\gamma} \left[ \id - P_\gamma(x) Q_\gamma(x) \right],
    \end{equation*}
    where we denote \(Q_\gamma(x)=I-\gamma \nabla^2 f(x)\), \(P_\gamma(x) := J \prox_{\gamma g}(x - \gamma \nabla f(x))\), and \(M_{\gamma}(x) \in \R^{n \times n}\) with \[
    \left[M_\gamma(x)\right]_{i,k} = \sum_{j = 1}^n \frac{\partial \left[Q_\gamma(x)\right]_{i,j}}{\partial_k} \left[R_\gamma\right]_j = - \gamma \sum_{j = 1}^n \frac{\partial \left[\nabla^2 f(x) \right]_{i, j}}{\partial_k} \left[ R_{\gamma} \right]_{j}.
    \]
\end{theorem}
\begin{proof}
    The claim follows from the product rule of differentiation applied to
    \begin{equation*}
        \nabla \fbe(x) = Q_\gamma(x) R_\gamma(x),
    \end{equation*}
    where the Jacobian of the fixed-point residual $JR_\gamma(x)$ immediately follows from application of the chain rule to
    \(
        R_\gamma(x) = \frac{1}{\gamma} \left[ x - \prox_{\gamma g}(x - \gamma \nabla f(x)) \right].
    \)
\end{proof}\noindent
At a critical point $x^\star$, the computation of $\nabla^2 \fbe(x^\star)$ only requires the second-order oracle of $f$ and $g$, i.e., \(\nabla^2 f\) and \(J\prox_{\gamma g}\), because $M_\gamma(x^\star) = 0$.
However, at points $x \in U_{x^\star}$ in a neighborhood around $x^\star$, the Hessian $\nabla^2 \fbe(x)$ also depends on $\nabla^3 f(x)$, quickly making its computation prohibitively expensive as the problem size increases. 
Another practical complication for designing second-order methods arises due to the fact that we have no way of knowing whether we are `sufficiently close' to a critical point for the Hessian to exist.  
To circumvent these difficulties, the following set-valued \emph{generalized Hessian}, initially proposed in \cite{themelis_acceleration_2019}, may prove useful
\begin{equation} \label{eq:generalized-hessian}
    \hat{\partial}^2 \fbe(x) \eqdef \left\{ \gamma^{-1} \mathcal{Q_\gamma} \left( \id - \mathcal{P_\gamma} \mathcal{Q_\gamma} \right) \bigg\vert ~ \substack{\mathcal{P_\gamma} \in \ \partial_C \prox_{\gamma g}(x - \gamma \nabla f(x))\\\mathcal{Q_\gamma} \in \partial_C(x - \gamma \nabla f(x))} \right\}.
\end{equation}
Remark that \eqref{eq:generalized-hessian} is well-defined and nonempty under \cref{assump:fbe-basic,assump:weak-convexity}, because weak convexity ensures local Lipschitz continuity of \(\prox_{\gamma g}\) for sufficiently small \(\gamma > 0\) \cite[Theorem 15.4.12]{themelis_acceleration_2019}.
This generalized Hessian is a Gauss-Newton-like approximant to $\nabla^2 \fbe$ which omits the term $M_\gamma$.
Since \(M_\gamma\) vanishes at critical points, $\hat{\partial}^2 \fbe$ reduces to (a singleton containing) the exact Hessian $\nabla^2 \fbe$ at critical points, as stated by the following corollary.
\begin{corollary} \label{corr:hessian-gn-connection}
    Suppose that \cref{assump:fbe-basic,assump:weak-convexity} hold.
    If \cref{assump:fbe-twice-diff} holds with respect to a critical point \( x^\star \in \R^n \) and a sufficiently small \(
    \gamma > 0\), then there exists a neighborhood $U_{x^\star}$ of $x^\star$ on which \(\hat{\partial}^2 \fbe(x)\) becomes a singleton \(\{B(x)\}\) where \(B(x) = \gamma^{-1} {Q_\gamma}(x) \left( \id - P_\gamma(x) {Q_\gamma}(x) \right)\).
    Moreover, for all \(x \in U_{x^\star}\), it holds that
    \begin{equation*}
        \Vert \nabla^2 \fbe(x) - B(x) \Vert = \Vert M_\gamma(x) \Vert \leq \gamma \Vert \nabla^3 f(x)\Vert \Vert R_\gamma(x) \Vert.
    \end{equation*}
\end{corollary}
A notable advantage of \eqref{eq:generalized-hessian} is that it does not involve third-order derivatives of \( f \), which is crucial for designing practical algorithms.
Observe also that the (generalized) Jacobian of the proximal mappings $\prox_{\gamma g}$ can in many cases be effectively computed.
We refer to \cite[\S 6]{themelis_acceleration_2019} for an extensive overview of closed-form expressions.

At this point, we have characterized the second-order differentiability of $\fbe$ around critical points in terms of the second-order epi-derivatives of $g$.
As a more practical alternative, we also show that the proximal operators of $\C^2$-partly smooth functions $g$ are locally continuously differentiable if a mild strict complementarity condition holds.
We emphasize that $\C^2$-partial smoothness, along with this strict complementarity condition, is a sufficient but not a necessary condition for local smoothness of the proximal operator.
Indeed, strict twice epi-differentiability in a neighborhood -- which is an equivalent characterization by \cref{th:prox-diff} -- does not amount to $\C^2$-partial smoothness \cite[Example 3.14]{hang_fresh_2024}.
\begin{theorem} \label{th:prox-c1-partly-smooth}
    Suppose that \cref{assump:fbe-basic} holds.
    Let $g$ be $\C^2$-partly smooth at a critical point $x^\star \in \R^n$ relative to a manifold \(\cM \subseteq \R^n\), and prox-regular at $x^\star$ for $-\nabla f(x^\star)$.
    If
    \begin{equation} \label{eq:strict-complementarity} \tag{SC}
        - \nabla f(x^\star) \in \relint{\partial g(x^\star)},
    \end{equation}
    then for all $\gamma > 0$ sufficiently small, the proximal mapping $\prox_{\gamma g}$ is continuously differentiable around the forward point $x^\star - \gamma \nabla f(x^\star)$.
\end{theorem}
\begin{proof}
    Consider the tilted function $\tilde{g} = g + \langle \nabla f(x^\star),\,\cdot\,\rangle$, which is prox-regular at $x^\star$ for $0$ \cite[Example 13.35]{rockafellar_variational_1998} and $\C^2$-partly smooth at $x^\star$ \cite[Corollary 4.6]{lewis_active_2002}.
    It follows from \cite[Theorem 28]{daniilidis_geometrical_2006} that the proximal operator $\prox_{\gamma \tilde{g}}$ is continuously differentiable around $x^\star$ for $\gamma$ sufficiently small.
    The claim follows by tilting back to $g$.
\end{proof}\noindent
Note that this theorem also follows directly from the recent result \cite[Corollary 3.19]{hang_fresh_2024}.

It is worth mentioning that trust-region methods for dealing with convex constraints require the strict complementarity condition \eqref{eq:strict-complementarity} to hold for each limit point to establish convergence to second-order stationary points, see \cite[AO.4b]{conn_trust_2000}.
Likewise, the second-order stationarity guarantees of the methods proposed in \cref{sec:nonsmooth-tr,sec:curvilinear} require \( \prox_{\gamma g} \) to be locally of class \( \C^1 \) around the so-called forward points of the limit points, which, in the case of \( \C^2 \)-partly smooth functions \( g \), reduces to the same assumption.

\begin{example}[Strict complementarity -- weakly active constraints] \label{ex:nonnegativity}
    Consider a problem of the form \eqref{eq:problem} where $g \equiv \delta_{\R^n_+}$ is the \textit{indicator function} of the \textit{nonnegative orthant} and $f$ is any sufficiently smooth function.
    The first-order optimality condition $-\nabla f(x^\star) \in \partial g(x^\star)$ thus becomes $-\nabla f(x^\star) \in N_{\R^n_+}(x^\star)$, and the strict complementarity condition \eqref{eq:strict-complementarity} becomes $- \nabla f(x^\star) \in \relint{N_{\R^n_+}(x^\star)}$.
    Here the normal cone $N_{\R^n_+}(x)$ and its relative interior $\relint{N_{\R^n_+}(x^\star)}$ are
    \begin{equation*}
        N_{\R^n_+}(x) = \left\{ v\ \middle| \ \begin{aligned} 
            v_i = 0 \quad x_i > 0\\
            v_i \leq 0 \quad x_i = 0
        \end{aligned} \right\}, \qquad \relint{N_{\R^n_+}(x)} = \left\{ v\ \middle| \ \begin{aligned} 
            v_i = 0 \quad x_i > 0\\
            v_i < 0 \quad x_i = 0
        \end{aligned} \right\},
    \end{equation*}
    i.e., any first-order stationary point $x^\star$ satisfies
    \begin{equation*}
        \begin{aligned}
            \left[ x^\star \right]_i = 0 \Rightarrow &\left[ \nabla f(x^\star) \right]_i \geq 0, \qquad
            \left[ x^\star \right]_i > 0 \Rightarrow &\left[ \nabla f(x^\star) \right]_i = 0.
        \end{aligned}
    \end{equation*}
    If the gradient components $\left[ \nabla f(x^\star) \right]_i > 0$ are strictly positive when $\left[ x^\star \right]_i = 0$, then also \eqref{eq:strict-complementarity} holds, and $\fbe$ is locally of class $\C^2$ around $x^\star$.
   Thus, \eqref{eq:strict-complementarity} states that there must not be any so-called \emph{weakly active} constraints, i.e., constraints for which \(\left[x^\star\right]_i = \left[ \nabla f(x^\star) \right]_i = 0\) \cite{nocedal_numerical_2006}.
\end{example}

%% file: fbe/equivalence-strict-saddles.tex
Having described generic conditions under which $\fbe$ is locally $\C^2$,
the question arises whether this result can be used to design methods that converge to \emph{second-order stationary points} of nonsmooth functions $\varphi$.
The following theorem shows that under \cref{assump:fbe-twice-epi}, which is the least restrictive variant of \cref{assump:fbe-twice-epi-all} and in particular a weaker condition than \cref{assump:fbe-twice-diff}, the second-order stationary points of $\varphi$ agree with those of $\fbe$.
\begin{theorem}[Equivalence of second-order stationary points] \label{th:equivalence-strict-saddles}
    Suppose that $\nabla^2 f$ exists and is continuous around a critical point $x^\star$, and that $g$ satisfies \cref{assump:fbe-twice-epi} at $x^\star$. Then, for $\gamma > 0$ sufficiently small we have that
    \begin{equation*}
        \forall s \in \bR^n: \mathrm{d}^2 \varphi(x^\star, 0)(s) \geq 0 \quad \Leftrightarrow \quad \forall s \in \bR^n: \langle \nabla^2 \fbe(x^\star) s, s \rangle \geq 0.
    \end{equation*}
    Consequently, the set of second-order stationary points of $\fbe$ and $\varphi$ are the same.
\end{theorem}
\begin{proof}
    By \cref{assump:fbe-twice-epi}, the second-order epi-derivative of $g$ is generalized quadratic, i.e.,
    \begin{equation} \label{eq:generalized-quadratic}
        \mathrm{d}^2 g(x^\star \mid-\nabla f(x^\star))[s]=\langle s, M s\rangle+\delta_S(s), \quad \forall s \in \mathbb{R}^n,
    \end{equation}
    where $S \subseteq \mathbb{R}^n$ is a linear subspace, and $M \in \mathbb{R}^{n \times n}$.
    Without loss of generality, we can take $M$ symmetric and such that $\operatorname{Im}(M) \subseteq S$ and $\operatorname{Ker}(M) \supseteq S^{\perp}$.\footnote{
        Remark that if $M$ and \(S\) satisfy \eqref{eq:generalized-quadratic}, then also $M' = \proj_S \left[ \frac{M + M^\top}{2} \right] \proj_S$ and \(S\) do.
        We can thus always replace $M$ by $M'$ to enforce the described properties.
    }
    Thus, for all $s \in \bR^n$ (see e.g.\,\cite[Eq. 6.1]{stella_forwardbackward_2017})
    \begin{align*}
        \mathrm{d}^2 \varphi(x^\star, 0)(s) &= \langle \nabla^2 f(x^\star) s, s \rangle + \langle M s, s \rangle + \delta_S(s),
    \end{align*}
    i.e.,
    \begin{equation*}
        \forall s \in \bR^n: \mathrm{d}^2 \varphi(x^\star, 0)(s) \geq 0 \quad \Leftrightarrow \quad \forall s \in S: \langle \nabla^2 f(x^\star) s, s \rangle + \langle M s, s \rangle \geq 0.
    \end{equation*}
    Moreover, invoking \cref{prop:fbe-twice-diff} and \(x^\star \in \fix T_\gamma\) ensures
    \begin{align*}
        \forall s \in \bR^n: \langle \nabla^2 \fbe(x^\star) s, s \rangle \geq 0 \quad &\Leftrightarrow \quad \forall s \in \bR^n: \langle \left[ Q_\gamma(x^\star) - Q_\gamma(x^\star) P_\gamma(x^\star) Q_\gamma(x^\star) \right] s, s \rangle \geq 0\\
        &\Leftrightarrow \quad \forall s' \in \bR^n: \langle Q_\gamma^{-1}(x^\star) s', s' \rangle - \langle P_\gamma(x^\star) s', s' \rangle \geq 0\\
        &\Leftrightarrow \quad \forall s' \in \bR^n:\langle Q_\gamma^{-1}(x^\star) s', s' \rangle - \langle \Pi_S (\id + \gamma M)^{-1} \Pi_S s', s' \rangle \geq 0.
    \end{align*}
    Here, we used the change of variables \(s' = Q_\gamma(x^\star) s\), and the fact that \(P_\gamma(x^\star) = \Pi_S (\id + \gamma M)^{-1} \Pi_S\) \cite[Lemma 2.9]{stella_forwardbackward_2017}.
    Observe that $Q_\gamma(x^\star)$ and $\id+\gamma M$ are both positive definite (see the intermediate result in the proof of \cite[Th 4.10i]{themelis_forward-backward_2018} for the latter).
    By using a Schur complement, we have that
    \(
        \nabla^2 \fbe(x^\star) \succeq 0 %
    \)
    if and only if \cite[Theorem 1.12]{zhang_schur_2005} (under nonsingularity of $\id+\gamma M$)
    \begin{equation*}
        \begin{bmatrix}
            \id+\gamma M & \Pi_S\\
            \Pi_S & Q_\gamma^{-1}(x^\star)
        \end{bmatrix} \succeq 0,
        \qquad \Leftrightarrow \qquad
        \begin{bmatrix}
            Q_\gamma^{-1}(x^\star) & \Pi_S\\
            \Pi_S & \id+\gamma M
        \end{bmatrix} \succeq 0,
    \end{equation*}
    if and only if \cite[Theorem 1.12]{zhang_schur_2005} (under nonsingularity of $Q_\gamma(x^\star)$)
    \begin{equation*}
        (\id + \gamma M) - \Pi_S Q_\gamma(x^\star) \Pi_S \succeq 0.  
    \end{equation*}
    Thus far, we have shown the equivalence
    \begin{equation*}
        \forall s \in \bR^n: \langle \nabla^2 \fbe(x^\star) s, s \rangle \geq 0 \Leftrightarrow (\id + \gamma M) - \Pi_S Q_\gamma(x^\star) \Pi_S \succeq 0.
    \end{equation*}
    We now first show the `$\Leftarrow$' part of the claim.
    Positive semidefiniteness of \((\id + \gamma M) - \Pi_S Q_\gamma(x^\star) \Pi_S\) on $\bR^n$ implies positive semidefiniteness on $S \subseteq \bR^n$ and thus one easily derives that
    \begin{equation*}
        \Pi_S (\id + \gamma M) \Pi_S - \Pi_S Q_\gamma(x^\star) \Pi_S \succeq 0,
    \end{equation*}
   equivalently 
    \begin{equation*}
        \Pi_S (M + \nabla^2 f(x^\star)) \Pi_S \succeq 0 \Leftrightarrow \forall s \in \bR^n: \mathrm{d}^2 \varphi(x^\star, 0)(s) \geq 0.
    \end{equation*}
    As for the direction `$\Rightarrow$', it holds that
    \begin{align*}
        &(\id + \gamma M) - \Pi_S Q_\gamma(x^\star) \Pi_S\\
        &=\Pi_S (\id + \gamma M) \Pi_S + \underbrace{\Pi_{S^\perp} (\id + \gamma M) \Pi_{S^\perp}}_{\succeq 0} + \Pi_S (\id + \gamma M) \Pi_{S^\perp} + \Pi_{S^\perp} (\id + \gamma M) \Pi_S - \Pi_S Q_\gamma(x^\star) \Pi_S\\
        &\geq \Pi_S (\id + \gamma M) \Pi_S + \gamma \left[ \left( \Pi_{S^\perp} M \Pi_S \right)^\top + \Pi_{S^\perp} M \Pi_S \right] - \Pi_S Q_\gamma(x^\star) \Pi_S\\
        &= \Pi_S (\id + \gamma M) \Pi_S - \Pi_S Q_\gamma(x^\star) \Pi_S \succeq 0.
    \end{align*}
    Here the last equality uses \(\Pi_{S^\perp} M \Pi_{S} = 0\), which follows from the fact that for any \(d \in \R^n\) we have \(M \Pi_S d \in S\) (due to \(\operatorname{Im}(M) \subseteq S\)) and hence \(\Pi_{S^\perp} M \Pi_{S} d = 0\).
\end{proof}\noindent
In light of \cref{th:equivalence-strict-saddles}, the original nonsmooth problem \eqref{eq:problem} is equivalent to
\begin{equation} \label{eq:problem-fbe} \tag{P-FBE}
    \minimize_{x \in \R^n} \fbe(x),
\end{equation}
in the sense that critical points coincide (cf.\,\cref{prop:fbe-decrease-minimizers}), strong local minimizers coincide (cf.\,\cite[Theorem 4.11]{themelis_forward-backward_2018}), and second-order stationary points coincide (cf.\,\cref{th:equivalence-strict-saddles}).
If \cref{assump:fbe-twice-diff} holds at critical points, then the FBE is locally of class $\C^2$, even though $\varphi$ need not be.
Thus, this equivalence provides a way for applying classical techniques to the smooth problem \eqref{eq:problem-fbe}, while guaranteeing convergence to a second-order stationary point of the original nonsmooth problem \eqref{eq:problem}.

We conclude by discussing some related results.
First, Liu and Yin \cite{liu_envelope_2019} establish a similar equivalence of second-order stationary points, but only under the more restrictive assumption that $g$ is locally $\C^2$.
Under this requirement, also the objective $\varphi$ is locally of class $\C^2$, a setting which, for example, does not cover any of the nonnegativity constraints in \cref{ex:nonnegativity} being active around the given critical point.
Second, Davis and Drusvyatskiy \cite{davis_proximal_2022} show that if \(\varphi\) admits a certain \emph{active manifold}\footnote{
    The active manifold in \cite{davis_proximal_2022} is different from that of \(\C^2\)-partly smooth functions:
    in certain settings its existence is equivalent to \(\C^2\)-partial smoothness under the constraint qualification \eqref{eq:strict-complementarity} \cite[Proposition 8.4]{drusvyatskiy_optimality_2014}.
} around a critical point, then the Moreau envelope of \(\varphi\) is locally of class \(\C^2\).
In addition, a second-order stationary point of \(\varphi\) -- which the authors define in terms of the active manifold -- is then also a second-order stationary point of its Moreau envelope.
We note that the analysis in \cite{davis_proximal_2022} heavily relies on this active manifold, whereas this work exploits the structure of the second-order epi-derivatives (cf.\,\cref{assump:fbe-twice-epi-all}).

%% file: fbtr/main.tex
Building on the second-order properties of the forward-backward envelope (FBE) established in the previous section, we now propose a nonsmooth trust-region method to solve problem \eqref{eq:problem}, equipped with second-order convergence guarantees akin to those of classical trust-region algorithms. Specifically, the method ensures convergence to second-order stationary points and achieves a locally superlinear rate of convergence.

\input{fbtr/algorithm.tex}

\input{fbtr/convergence.tex}

%% file: fbtr/algorithm.tex
\subsection{Algorithmic development}

A natural approach for designing a method that converges to second-order stationary points of \( \varphi \) is to leverage the equivalence between minimizing \( \varphi \) and minimizing the forward-backward envelope \( \fbe \), as established in \cref{sec:equivalence}.
In particular, if \( \fbe \) is sufficiently smooth, then classical second-order methods can be applied to \eqref{eq:problem-fbe} at least locally, and a second-order stationary point will also be a second-order stationary point of the original problem \eqref{eq:problem}.
Thus, this section investigates the idea of applying a classical \emph{trust-region method} to the surrogate problem \eqref{eq:problem-fbe}.
Trust-region methods\footnote{We refer to \cite{conn_trust_2000} for a thorough discussion of trust-region methods in the smooth setting.} iteratively minimize a local model of the objective in a neighborhood around the current iterate, typically based on a Taylor expansion.
To ensure sufficient smoothness of the FBE, \(f \in \C^{2+}(\R^n)\) is assumed throughout this section, along with \cref{assump:fbe-basic,assump:weak-convexity}.
\begin{assumption}\label{assump:fbe-basic-tr}
    \Cref{assump:fbe-basic,assump:weak-convexity} hold, and additionally $f \in \C^{2+}(\R^n)$.
\end{assumption}\noindent
\Cref{assump:fbe-basic-tr} ensures that $\fbe \in \C^{1+}$, which is a standard assumption used for ensuring global convergence of trust-region methods.
\begin{proposition}[Differentiability of the \(\fbe \)] \label{prop:fbe-C1}
    Suppose that \cref{assump:fbe-basic-tr} holds, and let $\gamma \in (0, \min \{ \nicefrac{1}{L_f}, \nicefrac{1}{\rho} \})$.
    Then, $\fbe \in \C^{1+}$ and the expression \( \nabla \fbe = Q_\gamma R_\gamma \) from \cref{prop:fbe-locally-c1} holds globally.
\end{proposition}
\begin{proof}
    Continuous differentiability of \( \fbe \) follows from \cite[Th.\,2.6]{stella_forwardbackward_2017}, while the local Lipschitz continuity of \( \nabla \fbe \) follows from \cite[Theorem 4.7(i)]{themelis_forward-backward_2018} and the fact that the product of locally Lipschitz continuous mappings is locally Lipschitz continuous.
\end{proof}\noindent
Let us begin by defining the quadratic subproblem
\begin{equation}\label{eq:subproblem} \tag{P-TR}
    \begin{aligned}
        &\minimize_{d} &&m_k(d) \eqdef \fbe(x^k) + \langle \nabla \fbe(x^k), d \rangle + \frac{1}{2} \langle B_k d, d \rangle\\
        &\stt &&\Vert d \Vert \leq \delta_k,
    \end{aligned}
\end{equation}
where $B_k \in \hat \partial^2 \fbe(x^k)$ and $\delta_k$ is the radius at iteration \(k \geq 0\).
This particular choice of Hessian approximant is fundamental to establishing convergence to second-order stationary points.
Conventional proofs require that, at least locally, the exact Hessian is used.
However, this approach would create multiple issues in our setting, namely: (i) $\nabla^2 \fbe$ does not exist everywhere; and (ii) even if $\nabla^2 \fbe$ exists locally, it is in general prohibitively expensive to compute since this involves third-order derivatives of $f$ (cf.\,\cref{prop:fbe-twice-diff}). 
By selecting $B_k \in \hat \partial^2 \fbe(x^k)$ as an element of the \emph{generalized Hessian} \eqref{eq:generalized-hessian}, both concerns are addressed, since \(B_k\) exists along the iterates and only involves second-order information of $f$ and $g$.
However, a difficulty in proving convergence to second-order stationary points in this way is that \(\hat \partial^2 \fbe(x^\star)= \{ \nabla^2 \fbe(x^\star) \}\) holds only at critical points \(x^\star \in \fix T_{\gamma}\), but not around them (cf.\,\cref{corr:hessian-gn-connection}).
The proof of \cref{thm:second-order-stationary} addresses this by exploiting the local structure of \(\nabla^2 \fbe\) and \(\hat \partial^2 \fbe\) around critical points.

Beyond the use of an unconventional Hessian approximant, we investigate the application of a standard trust-region method to the problem \eqref{eq:problem-fbe}.
An (approximate) minimizer $d^k$ of \eqref{eq:subproblem} constitutes a candidate update direction at the current iterate.
If the ratio $\trratio_k$ of the \textit{actual reduction} of $\fbe$ compared to the \textit{predicted reduction} of the model $m_k$,
\begin{equation} \label{eq:tr-ratio}
    \trratio_k = \frac{\fbe(x^k)-\fbe\left(x^k+d^k\right)}{m_k(0)-m_k\left(d^k\right)},
\end{equation}
is sufficiently large, i.e., if $\trratio_k \geq \mu_1$ for $\mu_1 > 0$, then the candidate direction $d^k$ is accepted.
The ratio $\trratio_k$ is also used to update the radius $\delta_k$.
We consider a simple update rule
\begin{equation} \label{eq:radius-rule}
    \begin{aligned}
        \delta_{k+1} = \left\{
            \begin{array}{ll}
                c_1 \delta_k & \quad \trratio_k < \mu_1,\\
                c_2 \delta_k & \quad \mu_1 \leq \trratio_k < \mu_2,\\
                c_3 \delta_k & \quad \trratio_k \geq \mu_2,
            \end{array}
        \right.
    \end{aligned}
\end{equation}
where $0<\mu_1<\mu_2<1$, $0 < c_1 < c_2 < 1 < c_3$.
The method is summarized in \cref{alg:fbtr}.

\begin{algorithm}
    \caption{NTRA (Nonsmooth Trust-Region Algorithm)}
    \label{alg:fbtr}
    \begin{algorithmic}[1]
        \Require{
            \indent {\small\(x_0\in\R^n\),
            \(\gamma\in(0,\min \left\{ \nicefrac{1}{L_f}, \nicefrac{1}{\rho} \right\})\), $0 < \delta_0$, $0<\mu_1<\mu_2<1$, $0 < c_1 < c_2 < 1 < c_3$};
        }
        \Initialize{
            $k \leftarrow 0$;
        }
        \For{$k=0, 1, \dots$}
            \State Select \( P_k \in \partial_C\left(\prox_{\gamma g})(x^k - \gamma \nabla f(x^k)\right) \);
            \State Let \( B_k = \gamma^{-1} Q_k (\id - P_k Q_k) \) with \( Q_k = \id - \gamma \nabla^2 f(x^k) \);
            \State \textbf{if} \( \Vert R_\gamma(x^k) \Vert = 0 \) and \( \lambda_{\min}(B_k) \geq 0 \) \textbf{then} stop; \textbf{end if}
            \State{Determine the direction $d^k$ as the (approximate) solution to \eqref{eq:subproblem};}
            \State{Compute the ratio $\trratio_k$ as in \eqref{eq:tr-ratio};}
            \If{$\trratio_k < \mu_1$}
                \State $x^{k+1} \leftarrow x^k$;
            \Else
                \State $x^{k+1} \leftarrow x^k + d^k$;
            \EndIf
            \State Update the trust-region radius $\delta_k$ according to \eqref{eq:radius-rule}.
        \EndFor
    \end{algorithmic}
\end{algorithm}\noindent

%% file: fbtr/convergence.tex
\subsection{Convergence analysis} \label{sec:convergence}

Smooth trust-region methods are typically analyzed under two key assumptions: (i) the objective function is level-bounded, and (ii) the sequence of Hessian approximants remains uniformly bounded; cf.\,\cite{nocedal_numerical_2006}.
The following two lemmas establish analogous results within our setting, under appropriate assumptions on the original problem \eqref{eq:problem}.
\begin{assumption} \label{ass:lowerLevelSet}
    The function $\varphi$ has bounded sublevel sets 
    \(
        \mathcal{L}(x^0):= \{x\in\R^n ~|~ \varphi(x)\leq \varphi(x^0)\}.
    \)
\end{assumption} \noindent
For the composite function $\varphi \equiv f + g$, \cref{ass:lowerLevelSet} is satisfied when one term is lower-bounded and the other term is level-coercive \cite[Exercise 3.29(a)]{rockafellar_variational_1998}.
\begin{lemma} \label{lem:bounded-level-sets-fbe}
    Suppose that \cref{assump:fbe-basic-tr,ass:lowerLevelSet} hold.
    Then \(\fbe\) with \(\gamma\in(0,\min \left\{ \nicefrac{1}{L_f}, \nicefrac{1}{\rho} \right\})\) has bounded sublevel sets $\mathcal{L}_\gamma(x^0) := \{ x \in \mathbb{R}^n ~|~ \fbe (x) \leq \fbe(x^0) \}$.
\end{lemma}
\begin{proof}
    By \eqref{eq:fbe-def} and \cref{prop:fbe-decrease-minimizers}, we have for any \(x \in \bR^n\) and \(\bar x \in T_\gamma(x)\) that
    \begin{equation*}
        f(x) + g(\bar x) + \langle \nabla f(x), \bar x - x \rangle + \frac{1}{2\gamma} \Vert \bar x - x \Vert^2 = \fbe(x).
    \end{equation*}
    Hence, the \(L_f\)-Lipschitz smoothness of \(f\) ensures
    \begin{align*}
        \frac{1}{2\gamma} \Vert \bar x - x \Vert^2 &\leq \fbe(x) - f(x) - g(\bar x) - \langle \nabla f(x), \bar x - x \rangle \leq \fbe(x) - \varphi(\bar x) + \frac{L_f}{2} \Vert \bar x - x \Vert^2.
    \end{align*}
    Moreover, the lower boundedness of \(\varphi\) yields 
    \begin{equation}\label{prf:bounded-level-sets-1}
        \frac{1 - \gamma L_f}{2\gamma} \Vert \bar x - x \Vert^2 \leq \fbe(x) - \inf \varphi.
    \end{equation}
    Now, by contradiction, suppose that there exists an unbounded sequence \(\seq{x^k}\) such that \(\Vert x^k \Vert \to \infty\) satisfying \(\fbe(x^k) \leq \alpha\) for some \(\alpha \geq \inf \varphi\) and all \(
    k \geq 0\).
    Then, by \eqref{prf:bounded-level-sets-1}, we have for all \(k \geq 0\)
    \begin{equation*}
        \Vert \bar x^k - x^k \Vert^2 \leq \frac{2\gamma}{1 - \gamma L_f} \left( \alpha - \inf \varphi \right) =: C.
    \end{equation*}
    Since \(\varphi(\bar x^k) \leq \fbe(x^k) \leq \alpha\) and since \(\varphi\) has bounded sublevel sets, there exists a \(R \geq 0\) such that \(\Vert \bar x^k \Vert \leq R\) for all \(k \geq 0\).
    By the triangle inequality, this yields
    \begin{equation*}
        \Vert x^k \Vert \leq \Vert \bar x^k - x^k \Vert + \Vert \bar x^k \Vert \leq \sqrt{C} + R,
    \end{equation*}
    which contradicts the assumption that \(\Vert x^k \Vert \to \infty\) and concludes the proof.
\end{proof}
\begin{lemma} \label{lem:bounded-Hessian-approximants} 
    Consider \(\seq{B_k}\) generated by \cref{alg:fbtr}.
    If \cref{assump:fbe-basic-tr} holds, then \(\Vert B_k \Vert \leq \nicefrac{6}{\gamma}\) holds for all \(k \geq 0\).
\end{lemma}
\begin{proof}
    We have \(\Vert B_k \Vert = \gamma^{-1} \Vert Q_k - Q_k P_k Q_k \Vert \leq \gamma^{-1} \Vert Q_k \Vert + \gamma^{-1} \Vert Q_k \Vert^2 \Vert P_k \Vert\).
    By Lipschitz-smoothness of \(f\) it follows for all \(x \in \R^n\) that \(\Vert \nabla^2 f(x) \Vert \leq L_f\) and hence \(\Vert I - \gamma \nabla^2 f(x) \Vert \leq 1 + \gamma L_f < 2.\)
    We conclude that \(\Vert Q_k \Vert \leq 2\) for \(k \geq 0\).
    On the other hand, \cite[Theorem 15.4.12]{themelis_acceleration_2019} immediately implies that \(\Vert P_k \Vert \leq 1\) for \(k \geq 0\), from which the claim follows.
\end{proof}\noindent
It is common to assume a sufficient decrease of the model $m_k$ by the (approximate) solution $d^k$ of \eqref{eq:subproblem}; see, e.g., \cite{nocedal_numerical_2006}.
The following sufficient decrease condition holds for various popular inner solvers, including the Dogleg method and Steihaug's conjugate gradient method \cite{nocedal_numerical_2006}.
\begin{condition}[Sufficient decrease of the model]\label{assump:suff-decrease}
    We assume that the decrease on the model $m_k$ is at least a fraction of the decrease obtained by the Cauchy point, i.e.,
    \begin{equation*}
        m_k(0) - m_k(d^k) \geq \beta \Vert \nabla \fbe(x^k) \Vert \min \left\{ \delta_k, \frac{\Vert \nabla \fbe(x^k) \Vert}{\Vert B_k \Vert} \right\}
    \end{equation*}
    for some constant $\beta \in (0, 1)$ and for all $k$.
\end{condition}\noindent
Given \cref{lem:bounded-level-sets-fbe,lem:bounded-Hessian-approximants}, well-definedness and global convergence of \cref{alg:fbtr} follow by standard results. 
For completeness, a proof is provided in \cref{app:well-defined-global}.
\begin{theorem}[Global convergence] \label{th:fbe-limit}
    Suppose that \cref{assump:fbe-basic-tr,ass:lowerLevelSet} and \cref{assump:suff-decrease} hold. 
    Then, $\left( x^k \right)_{k \in \N}$ generated by \cref{alg:fbtr} satisfies
    \begin{equation*}
        \lim_{k \to \infty} \Vert \nabla \fbe (x^k) \Vert = 0.
    \end{equation*}
\end{theorem}

To guarantee convergence to second-order stationary points, we additionally require a condition ensuring that the approximate solution $d^k$ to the subproblem \eqref{eq:subproblem} captures sufficient negative curvature, if any exists.
As such, we impose a similar condition as in \cite[Condition 2]{shultz_family_1985}.

\begin{condition} \label{assumption:shultz-condition-two}
    If $B_k$ is indefinite, then the step $d^k$ yields a model decrease that is at least as good as a direction of sufficient negative curvature, i.e.
    \begin{equation*}
        m_k(0) - m_k(d^k) \geq \beta_2 (-\lambda_{min}(B_k)) \delta_k^2
    \end{equation*}
    for some constant $\beta_2 > 0$ and all $\delta_k > 0$.
\end{condition}\noindent
The subsequent result formally states that \cref{alg:fbtr} converges to second-order stationary points of \cref{eq:problem}.
We note that the proof follows a strategy similar to classical trust-region methods but additionally accounts for the fact that, in general, \( B_k \in \hat \partial^2 \fbe(x^k) \) coincides with the exact Hessian \( \nabla^2 \fbe(x^k) \) only at critical points. 
\begin{theorem}[Convergence to second-order stationary points] \label{thm:second-order-stationary}
    Suppose that the iterates generated by \cref{alg:fbtr} converge to a critical point $x^\star \in \R^n$, concerning which \cref{assump:fbe-twice-diff} holds, and that \cref{assump:fbe-basic-tr,ass:lowerLevelSet} hold.
    If the iterates satisfy \cref{assumption:shultz-condition-two}, 
    then $x^\star$ is a \emph{second-order stationary point} of \eqref{eq:problem}.
\end{theorem}

\begin{proof}
    By \cref{th:equivalence-strict-saddles}, $x^\star$ is a second-order stationary point of \eqref{eq:problem} if and only if it is a second-order stationary points of \eqref{eq:problem-fbe}, i.e., $H_\star := \nabla^2 \fbe(x^\star) \succeq 0$.
    Henceforth, we assume that $k$ is large enough, such that $H_k \eqdef \nabla^2 \fbe(x^k)$ is well-defined and continuous.
    It follows from the mean-value theorem that
    \begin{equation*}
            \fbe(x^k + d^k) = \fbe(x^k) + \langle \nabla \fbe(x^k), d^k \rangle + \frac{1}{2} \langle H(x^k) d^k, d^k \rangle + \frac{1}{2} \int_{0}^{1} \langle (H(x^k + t d^k) - H(x^k)) d^k, d^k \rangle dt,
    \end{equation*}
    i.e.,
    \begin{equation} \label{prf:fbe-second-order-taylor-bound}
        \begin{aligned}
            \vert \fbe(x^{k}) &- \fbe(x^{k} + d^{k}) - (m_{k}(0_n) - m_{k}(d^{k}))\vert\\
            &= \bigg\vert \frac{1}{2} \langle (B_k - H(x^k)) d^k, d^k \rangle - \frac{1}{2} \int_{0}^{1} \langle (H(x^k + t d^k) - H(x^k)) d^k, d^k \rangle dt \bigg\vert\\
            &\leq \frac{\gamma}{2} \Vert \nabla^3 f(x^k) \Vert \cdot \Vert R_\gamma(x^k) \Vert \cdot \Vert d^k \Vert^2 + o(\Vert d^k \Vert^2).
        \end{aligned}
    \end{equation}
    \emph{Suppose for contradiction,} that $\lambda_\star := \lambda_{\min}(H_\star) < 0$.
    Since by \cref{prop:fbe-twice-diff}, $B_k \to H_\star$, there exists an index $K > 0$ such that for $k \geq K$ we have
    \(
        \lambda_\text{min}(B_k) \leq \frac{\lambda_\star}{2} < 0.
    \)
    Under \cref{assumption:shultz-condition-two}, it then holds for $k \geq K$ that
    \begin{equation} \label{prf:fbe-second-order-prediction-bound}
        m_k(0_n) - m_k(d^k) \geq \beta_2 (-\lambda_\text{min}(B_k)) \delta_k^2 \geq \frac{\beta_2}{2} (-\lambda_\star) \delta_k^2.
    \end{equation}
    Thus, for $k \geq K$, it follows from \eqref{prf:fbe-second-order-taylor-bound} and \eqref{prf:fbe-second-order-prediction-bound} that 
    \begin{equation*}
        \begin{aligned}
            \lvert \trratio_{k} - 1 \rvert &\overset{\eqref{eq:tr-ratio}}{=} \bigg\lvert \frac{\fbe(x^{k}) - \fbe(x^{k} + d^{k})}{m_{k}(0_n) - m_{k}(d^{k})} - 1\bigg\rvert = \bigg\lvert \frac{\fbe(x^{k}) - \fbe(x^{k} + d^{k}) - (m_{k}(0_n) - m_{k}(d^{k}))}{m_{k}(0_n) - m_{k}(d^{k})}\bigg\rvert\\
            &\leq \frac{\frac{\gamma}{2} \Vert \nabla^3 f(x^k) \Vert \cdot \Vert R_\gamma(x^k) \Vert \cdot \Vert d^k \Vert^2 + o(\delta_k^2)}{\frac{\beta_2}{2} (-\lambda_\star) \delta_k^2} \leq \frac{\gamma \Vert \nabla^3 f(x^k) \Vert \cdot \Vert R_\gamma(x^k) \Vert}{\beta_2 (-\lambda_\star)} + \frac{o(\delta_k^2)}{\delta_k^2}.
        \end{aligned}
    \end{equation*}
    We now follow a reasoning similar to that in the proof of \cite[Theorem 2.2]{shultz_family_1985}.
    The first term on the right-hand side vanishes as $k \to \infty$, by criticality of $x^\star$.
    Moreover, the second term goes to $0$ as $\delta_k \to 0$.
    We conclude that the radius $\delta_k$ must be bounded below, since for sufficiently small $\delta_k$ the above inequality implies $\trratio_k \geq \mu_2$, which in turn prevents further decrease of \(\delta_k\) -- we have $\delta_{k+1} \geq \delta_k$ by the update rule \eqref{eq:radius-rule}.
    Lower boundedness of $\delta_k$ guarantees the existence of an infinite subsequence $\cK \subseteq \N$ such that $\trratio_k \geq \mu_1, \forall k \in \cK$.
    Then, it follows from \eqref{prf:fbe-second-order-prediction-bound} that for $k \in \mathcal{K}$ and $k \geq K$
    \begin{equation*}
            \fbe(x^k) - \fbe(x^{k+1}) = \fbe(x^k) - \fbe(x^k + d^k) \overset{\trratio_k \geq \mu_1}{\geq} \mu_1 (m_k(0_n) - m_k(d^k)) \geq \mu_1 \frac{\beta_2}{2} (-\lambda_\star) \delta_k^2.
    \end{equation*}
    Since $\fbe$ is bounded below, a standard telescoping argument yields
    \(
        \lim_{\substack{k \to \infty\\ k \in \mathcal{K}}} \delta_k = 0,
    \)
    contradicting the lower boundedness of $\delta_k$.
    This means that the original assertion $\lambda_\star < 0$ is false, giving our desired results.
\end{proof}

A local superlinear rate of convergence can be achieved if the inner solver computes the exact Newton step whenever it lies within the radius \( \delta_k \).
As before, this proof follows a strategy similar to classical trust-region methods but also addresses the fact that \( B_k \in \hat \partial^2 \fbe(x^k) \) coincides with the exact Hessian \( \nabla^2 \fbe(x^k) \) only at critical points.. 
\begin{condition}[{\cite[Condition 3]{shultz_family_1985}}] \label{assumption:shultz-condition-three}
    If the approximate Hessian $B_k$ is positive definite and the Newton step \(d_N^k := - B_k^{-1} \nabla \fbe(x^k) \) is sufficiently small, i.e.\,$\Vert d_N^k \Vert \leq \delta_k$, then the inner solver returns this Newton step as a solution, meaning that
    \begin{equation*}
        d^k = d_N^k.%
    \end{equation*}  
\end{condition}

\begin{theorem}[Local superlinear convergence rate] \label{th:quadratic-convergence}
    Suppose that that the iterates generated by \cref{alg:fbtr} converge to a critical point $x^\star \in \R^n$, with respect to which \cref{assump:fbe-twice-diff} holds, and that \cref{assump:fbe-basic-tr,ass:lowerLevelSet} hold.
    If $x^\star$ is a strong local minimizer of $\varphi$, and if \cref{assumption:shultz-condition-three} holds,
    then 
    \begin{theoremenum}
        \item the full sequence of iterates $\seq{x^k}$ converges \emph{sequentially} to $x^\star$;
        \item the Newton step is eventually accepted, and the iterates converge with a superlinear rate.
    \end{theoremenum}
\end{theorem}
\begin{proof}
    We denote the (unconstrained) Newton step by $d_N^k = -B_k^{-1} \nabla \fbe(x^k)$, whereas \(d^k\) indicates the actual step taken by \cref{alg:fbtr}.
    By \cref{corr:hessian-gn-connection}, for $x$ sufficiently close to $x^\star$ the generalized Hessian is a singleton \( \hat \partial^2 \fbe(x) = \left\{ B(x) \right\} \), where \[
        B(x) := \gamma^{-1} Q_\gamma(x) (\id - P_\gamma(x) Q_\gamma(x)), \qquad P_\gamma(x) := J \prox_{\gamma g}(x - \gamma \nabla f(x)),
    \]
    with \( B(x) \) continuous. For \( x^k \) sufficiently close to \( x^\star \) we have that \( B_k = B(x^k) \).\\

    \noindent
    Sequential convergence to $x^\star$ follows by similar arguments as in \cite[Theorem 2.2]{shultz_family_1985}.
    We first introduce two radii $\delta_1, \delta_2 >0$.
    Since $B_\star := \nabla^2 \fbe(x^\star)$ is positive definite and $B(x)$ is continuous in a neighborhood around $x^\star$, there exists a $\delta_1 > 0$ such that \textit{(i)} $\Vert x^k - x^\star \Vert < \delta_1$ results in $B_k \succ 0$ and \textit{(ii)} $x^k \neq x^\star$ results in $\nabla \fbe(x^k) \neq 0$.
    Denote $\cN_1 \eqdef \left\{ x \mid \Vert x - x^\star \Vert < \delta_1 \right\}$ and observe that $x^\star$ is the unique minimizer of $\fbe$ in $\cN_1$.
    Since $\nabla \fbe(x^\star) = 0$, we can additionally define a $\delta_2 > 0$ such that \textit{(i)} $\delta_2 < \nicefrac{\delta_1}{4}$ and \textit{(ii)} $\Vert B_k^{-1} \nabla \fbe(x^k) \Vert < \nicefrac{\delta_1}{2}$ for all $x \in \cN_2 \eqdef \left\{ x \mid \Vert x - x^\star \Vert < \delta_2 \right\}$.
    Note that $x^\star$, the unique minimizer of $\fbe$ in $\cN_1$, lies in $\cN_2$ too.
    The subsequential convergence of the iterates to $x^\star$ implies that there exists a $k_0 \in \N$ such that $\fbe(x^{k_0}) < \inf \left\{ \fbe(x) \mid x \in \cN_1 \setminus \cN_2 \right\}$ with $x^{k_0} \in \cN_2$.
    We now claim that for any $x^k$ with $k \geq k_0$ and $x^k \in \cN_2$, also $x^{k+1} \in \cN_2$, which implies that the entire sequence of iterates beyond $x^{k_0}$ lies in $\cN_2$.
    Suppose \emph{by contradiction} that $x^{k+1} \notin \cN_2$, then since $\fbe(x^{k+1}) < \fbe(x^{k_0})$ (due to the monotone updates of \cref{alg:fbtr}), $x^{k+1}$ is not in $\cN_1$, either, yielding,
    \begin{equation*}
        \delta_k \geq \Vert x^{k+1} - x^k \Vert \geq \Vert x^{k+1} - x^\star \Vert - \Vert x^k - x^\star \Vert \geq \delta_1 - \frac{\delta_1}{4} = \frac{3}{4} \delta_1 > \frac{\delta_1}{2} > \Vert B_k^{-1} \nabla \fbe(x^k) \Vert.
    \end{equation*} 
    \cref{assumption:shultz-condition-three} ensures that $d^k = - B_k^{-1} \nabla \fbe(x^k)$ and therefore $\Vert d_k \Vert < \nicefrac{\delta_1}{2}$.
    Combined with the fact that \(\Vert x^k - x^\star \Vert < \frac{\delta_1}{4}\), which follows from \(x^k \in \cN_{2}\), we obtain
    \begin{equation*}
        \Vert x^{k+1} - x^\star \Vert \leq \Vert x^k - x^\star \Vert + \Vert d_k \Vert < \frac{3}{4} \delta_1,
    \end{equation*}
    which implies $x^{k+1} \in \cN_1$, and results in a contradiction.
    Thus, $x^k \in \cN_2$ for all $k \geq k_0$.
    Since $\fbe(x^k)$ decreases monotonically along the iterates and a subsequence of the iterates converges to the unique minimizer of $\cN_2$, i.e., $x^\star$, it follows that the entire sequence $\seq{x^k}$ converges to $x^\star$.\\

    \noindent
    We now show that eventually the Newton step lies within the radius and is accepted.
    For $k$ sufficiently large, $B_k$ is positive definite, so that by \cref{assumption:shultz-condition-three}, either $\Vert d_N^k \Vert > \delta_k$, or $d^k = d_N^k$.
    In either case, $\Vert d^k \Vert \leq \Vert d_N^k \Vert \leq \Vert B_k^{-1} \Vert \Vert \nabla \fbe(x^k) \Vert$, and hence $\Vert \nabla \fbe(x^k) \Vert \geq \Vert d^k \Vert / \Vert B_k^{-1} \Vert$.
    From \cref{assump:suff-decrease} we then obtain that for all sufficiently large iterates \(k\)
    \begin{equation*}
        \begin{aligned}
            m_k(0_n) - m_k(d^k) &\geq \beta \Vert \fbe(x^k) \Vert \min \left\{ \delta_k, \frac{\Vert \nabla \fbe(x^k) \Vert}{\Vert B_k \Vert} \right\} 
            \geq \beta \frac{\Vert d^k \Vert}{\Vert B_k^{-1} \Vert} \min \left\{ \Vert d^k \Vert, \frac{\Vert d^k \Vert}{\Vert B_k \Vert \Vert B_k^{-1} \Vert} \right\}\\
            &= \beta \frac{\Vert d^k \Vert^2}{\Vert B_k^{-1} \Vert^2 \Vert B_k \Vert}
            \geq \beta \frac{\Vert d^k \Vert^2}{2 \Vert B_\star^{-1} \Vert^2 \Vert B_\star \Vert} = C \Vert d^k \Vert^2,
        \end{aligned}
    \end{equation*}
    where $C = \frac{\beta}{2 \Vert B_\star^{-1} \Vert^2 \Vert B_\star \Vert}$.
    Here, the last inequality used $\Vert B_k^{-1} \Vert^2 \Vert B_k \Vert \leq 2 \Vert B_\star^{-1} \Vert \Vert B_\star \Vert$ for all $k$ large enough.
    The same arguments that led to \eqref{prf:fbe-second-order-taylor-bound} still hold, so following a similar reasoning as in the proof of \cref{thm:second-order-stationary}, we find that \(\lvert \trratio_k - 1 \rvert \to 0 \) as \( \delta_k \to 0 \), from which 
    we conclude that $\delta_k$ is lower bounded.
    Since $d_N^k \to 0$, the Newton directions satisfy $\Vert d_N^k \Vert \leq \delta_k$ for $k$ large enough, and hence by \cref{assumption:shultz-condition-three} we find that eventually $d^k = d_N^k$.
    Observing \( m_k(0_n) - m_k(d_N^k) = \frac{1}{2} \langle B_k d_N^k, d_N^k \rangle \geq \lambda_{\min}(B_\star) \Vert d_N^k \Vert^2 \) and again using \eqref{prf:fbe-second-order-taylor-bound}, we find $\trratio_k \to 1$ as $\Vert d_N^k \Vert \to 0$.
    Thus, eventually $d^k = d_N^k$ is also accepted by \cref{alg:fbtr}.\\

    \noindent
    Finally, we show the superlinear rate.
    It follows from $d_N^k = -B_k^{-1} \nabla \fbe(x^k)$ that
    \begin{equation*}
        x^k + d_N^k - x^\star = B_k^{-1} \left[ B_k (x^k - x^\star) - (\nabla \fbe(x^k) - \nabla \fbe(x^\star)) \right].
    \end{equation*}
    Denote \( H_k := \nabla^2 \fbe(x^k) \) for \( k \) sufficiently large.
    Invoking the mean-value theorem ensures
    \begin{equation*}
        \nabla \fbe(x^k) = \nabla \fbe(x^\star) + H_k (x^k - x^\star) + \int_{0}^{1} [H(x^\star + t(x^k - x^\star)) - H_\star](x^k - x^\star) \mathrm d t.
    \end{equation*}
    By continuity of \(H(x) := \nabla^2 \fbe(x)\) and the inequality
    \begin{equation*}
        \bigg \Vert \int_{0}^{1} [H(x^\star + t(x^k - x^\star)) - H_\star](x^k - x^\star) \mathrm d t \bigg \Vert \leq \int_{0}^{1} \Vert H(x^\star + t(x^k - x^\star)) - H_\star \Vert \Vert x^k - x^\star \Vert \mathrm d t,
    \end{equation*}
    it is clear that
    \(
        \Vert H_k (x^k - x^\star) - \left( \nabla \fbe(x^k) - \nabla \fbe(x^\star) \right) \Vert = o (\Vert x^k - x^\star \Vert).
    \)
    Moreover, it follows by \cref{corr:hessian-gn-connection} that for points $x$ sufficiently close to $x^\star$
    \begin{equation*}
        \Vert (B_k - H_k) (x^k - x^\star) \Vert = \Vert M_\gamma(x^k) (x^k - x^\star) \Vert \leq \gamma \Vert \nabla^3 f(x^k) \Vert \Vert R_\gamma(x^k) \Vert \Vert x^k - x^\star \Vert.
    \end{equation*}
    By Lipschitz continuity of $R_\gamma(x)$ around $x^\star$, we find
    \(
        \Vert (B_k - H_k) (x^k - x^\star) \Vert = O ( \Vert x^k - x^\star \Vert^2 ).
    \)
    Nonsingularity of $B_\star$ implies that $\Vert B_k^{-1} \Vert \leq 2 \Vert B_\star^{-1} \Vert$ for all $k$ large enough.
    In conclusion, we find that
    \begin{equation*}
        \begin{aligned}
            \Vert x^k + d_N^k - x^\star \Vert &\leq \Vert B_k^{-1} \Vert \Vert B_k (x^k - x^\star) - (\nabla \fbe(x^k) - \nabla \fbe(x^\star)) \Vert\\
            &= \Vert B_k^{-1} \Vert \Vert [ ( B_k - H_k ) (x^k - x^\star) ] + [ H_k(x^k - x^\star) - ( \nabla \fbe(x^k) - \nabla \fbe(x^\star) ) ] \Vert\\
            &= o(\Vert x^k - x^\star \Vert),
        \end{aligned}
    \end{equation*}
    leading to the local superlinear convergence rate.
\end{proof}

\begin{remark}[Local quadratic convergence]

    We highlight that if $\nabla^2 \fbe$ is locally Lipschitz continuous around $x^\star$\footnote{
        Local Lipschitz continuity of $P_\gamma$ and $\nabla^3 f$ around $x^\star$ would guarantee this.
    }, then \cref{alg:fbtr} actually obtains a local \emph{quadratic} rate of convergence. 
    Indeed, this additional assumption ensures that \[
        \Vert H_k (x^k - x^\star) - \left( \nabla \fbe(x^k) - \nabla \fbe(x^\star) \right) \Vert = O (\Vert x^k - x^\star \Vert^2),
    \]
    and following the proof of \cref{th:quadratic-convergence} this yields
    \(
        \Vert x^k + d_N^k - x^\star \Vert = O(\Vert x^k - x^\star \Vert^2).
    \)
    This result is in line with the local quadratic convergence rate of Gauss-Newton methods whenever the approximant becomes exact at the (stationary) limit point.
\end{remark}

%% file: curvilinear/main.tex
This section introduces a curvilinear linesearch method that converges to second-order stationary points of \cref{eq:problem} and achieves a local superlinear convergence rate.
Unlike the trust-region method presented in \cref{sec:nonsmooth-tr}, this algorithm is not derived from the equivalent problem \eqref{eq:problem-fbe} but instead directly addresses \eqref{eq:problem}. Moreover, it applies to a broader class of problems since, in particular, \cref{ass:lowerLevelSet} is not required to guarantee global convergence.

\subsection{Algorithmic development}

Beyond trust-region methods, there is also a class of globally convergent linesearch algorithms that reliably converge to second-order stationary points.
Following the studies in \cite{mccormick_modification_1977,more_use_1979}, {\it curvilinear linesearch methods} employ a pair of descent directions $(d, s)$, where, broadly speaking, $d$ is a descent direction derived from positive curvature information of the Hessian, and $s$ corresponds to a direction of negative curvature.
The linesearch procedure then searches along a \emph{curvilinear} path, typically chosen of the form \[
    x(\tau) = x + \tau^2 d + \tau s.
\]
Traditional variants require the objective function to be twice continuously differentiable $\C^2$, a property that the forward-backward envelope (FBE) satisfies only locally. Consequently, existing schemes cannot be applied directly to \eqref{eq:problem-fbe}.
Instead, this section generalizes the \zerofpr{} method \cite{themelis_forward-backward_2018} for solving \eqref{eq:problem} to incorporate a curvilinear linesearch procedure, and establishes second-order convergence guarantees.
The proposed method searches along a curvilinear path \[
x(\tau) = \bar x + \tau^2 d + \tau s, \qquad \bar x \in T_\gamma(x).
\]
When using the FBE as a merit function, this approach ensures global convergence under \cref{assump:fbe-basic,assump:weak-convexity}, since \( \fbe \) is strictly continuous and satisfies \(
    \fbe(\bar x) \leq \fbe(x) - \sigma \Vert r \Vert^2\) , where \( r \in R_\gamma(x) 
\)
(cf.\,\cref{prop:fbe-decrease-minimizers}).

The particular choice of directions \( d\) and \(s\) determines whether convergence to second-order stationary points and a local superlinear convergence rate is obtained.
Classical curvilinear linesearch methods require \( d \) to be a descent direction\footnote{
    Classical variants additionally require \( d \) to be a \emph{strict} descent direction whenever the Hessian has no negative eigenvalues. 
    However, the proposed method does not need such a condition, because a sufficient decrease can be guaranteed for the curvilinear path involving \( \bar x \in T_\gamma(x) \). 
}, and likewise we require
\begin{equation*}
    \langle Q \bar r, d \rangle \leq 0, \qquad \text{for some} \qquad Q \in \id - \gamma \partial_C(\nabla f)(\bar x) \qquad \text{and} \qquad \bar r \in R_\gamma(\bar x). 
\end{equation*}
Whenever \( \fbe \) is differentiable at \( \bar x \), this reads \( \langle \nabla \fbe(\bar x), d \rangle \leq 0 \), meaning that \(d\) is a descent direction of the merit function \(\fbe\) at \( \bar x \).
We require the direction \( s \) to satisfy a similar descent criterion, and, additionally, it must be a direction of negative curvature of \( B \in \hat \partial^2 \fbe(\bar x) \), i.e.,
\begin{equation*}
    \langle Q \bar r, s \rangle \leq 0, \qquad \text{and} \qquad \begin{cases}
        \langle B s, s \rangle < 0 & \text{if } \lambda_{\min}(B) < 0,\\
        \langle B s, s \rangle = 0 & \text{otherwise}.
    \end{cases}
\end{equation*}

The proposed linesearch condition combines the one from \zerofpr{}, which is based on the decrease of \( \fbe \) obtained by a forward-backward step, and the one from classical curvilinear linesearch methods, which is based on a second-order Taylor expansion.
In particular, given parameters \( \mu \in (0, 1) \) and \( \sigma \in (0, \gamma \frac{1 - \gamma L_f}{2}) \) the linesearch procedure searches the largest nonnegative \( \tau \leq 1 \) for which
\begin{equation*}
    \fbe(x(\tau)) \leq \fbe(x) - \sigma \Vert r \Vert^2 + \frac{\mu}{2} \tau^2 \langle B s, s \rangle, \qquad \text{where} \qquad r \in R_\gamma(x).
\end{equation*}
The resulting algorithm is summarized in \cref{alg:curvilinear-panoc}.
\begin{algorithm}
    \caption{PGCL (Proximal Gradient algorithm with Curvilinear Linesearch)}
    \label{alg:curvilinear-panoc}
    \begin{algorithmic}[1]
        \Require{
            \(x_0\in\R^n\),
            \(\gamma\in(0,\min \left\{ \nicefrac{1}{L_f}, \gamma_g \right\})\), $\sigma \in (0, \gamma \frac{1 - \gamma L_f}{2})$, $\beta \in (0, 1)$, $\mu \in (0, 1)$.
        }
        \For{$k=0, 1, \dots$}
            \State{
                Select $\bar x^k \in T_\gamma(x^k)$, and set $r^k = \frac{1}{\gamma}(x^k - \bar x^k)$;
            }
            \State{
                Select $\bar r^k \in R_\gamma(\bar x^k)$, \( Q_k \in \id - \gamma \partial_C (\nabla f)(\bar x^k) \) and \( P_k \in \partial_C (\prox_{\gamma g})(\bar x^k - \gamma \nabla f(\bar x^k)) \);
            }
            \State{
                Define $B_k = \gamma^{-1} Q_k (\id - P_k Q_k)$; %
            }
            \State{
                \textbf{if} $\Vert \bar r^k \Vert = 0$ and $\lambda_{\text{min}}(B_k) \geq 0$ \textbf{then} stop; \textbf{end if}
            }
            \State{
                Select directions $d^k, s^k \in \R^n$ satisfying
                \begin{equation*} \label{eq:negative-curvature-direction}
                    \langle Q_k \bar r^k, d^k \rangle \leq 0, \qquad \langle Q_k \bar r^k, s^k \rangle \leq 0, \qquad \begin{cases}
                        \langle B_k s^k, s^k \rangle < 0 & \text{if } \lambda_{\min}(B_k) < 0\\
                        \langle B_k s^k, s^k \rangle = 0 & \text{otherwise};
                    \end{cases}
                \end{equation*}
            }
            \State{
                Let $\tau_k \in \left\{ \beta^m \mid m \in \N \right\}$ be the largest such that 
                $x^{k+1} = \bar x^k + \tau_k^2 d^k + \tau_k s^k$ satisfies
                \begin{equation} \label{eq:curvilinear-linesearch}
                    \fbe(x^{k+1}) \leq \fbe(x^k) - \sigma \Vert r^k \Vert^2 + \frac{\mu}{2} \tau_k^2 \langle B_k s^k, s^k \rangle.
                \end{equation}
            }
        \EndFor
    \end{algorithmic}
\end{algorithm}\noindent

A simple choice for the negative curvature direction $s^k$ is the following:
if $\lambda_{\text{min}}(B_k) < 0$, select $s^k = \pm v^k$ where $v^k$ is an eigenvector corresponding to the most negative eigenvalue of $B_k$, with the sign such that $\langle Q_k \bar r^k, s^k \rangle \leq 0$;
otherwise, set $s^k = 0$.
This eigenvector \( v^k \) can be computed using the Lanczos algorithm, which only requires matrix-vector products with \( B_k \).
Besides being memory-friendly, an additional advantage of the Lanczos algorithm is that it is intimately related to the conjugate gradient method, and therefore provides an appealing option for computing the positive curvature direction \( d^k \).
We refer to \cite[\S4]{gould_exploiting_2000} and the references therein for more details.

Finally, it appears that the relative scaling of the directions \( d^k \) and \( s^k \) has a significant impact, both on the overall performance of the method, and on the quality of the solutions (cf.\,\cref{table:phase_retrieval_best_global} below).
For this reason, when \( \lambda_{\min}(B_k) < 0 \), we use negative curvature directions \( s^k = \pm \rho_k v^k \), where \( v^k \) is an (approximate) eigenvector corresponding to the most negative eigenvalue of \( B_k \), and where \( \rho_k > 0 \) is a suitable scaling factor, defined as
\begin{equation*}
    \rho_k = \bar s \sqrt{-\lambda_{\min}(B_k)} \min \left\{ 1, \frac{1}{\Vert Q_k \bar r^k \Vert} \right\},
\end{equation*}
for some fixed \( \bar s > 0 \).
For \( \bar s = 1 \), which is used unless mentioned otherwise, we recover the scaling of \cite{lucidi_curvilinear_1998}.

\subsection{Convergence analysis}

Before analyzing the convergence properties of \cref{alg:curvilinear-panoc}, the following lemma establishes that the linesearch procedure terminates after a finite number of backtracking steps.
\begin{lemma}[Well-definedness of the linesearch]
    Suppose that \cref{assump:fbe-basic,assump:weak-convexity} hold.
    Then, for all \( k \in \N \) there exists a \( \bar \tau_k > 0 \) such that
    \begin{equation*}
        \fbe(\bar x^{k} + \tau^2 d^k + \tau s^k) \leq \fbe(x^k) - \sigma \Vert r^k \Vert^2 + \frac{\mu}{2} \tau^2 \langle B_k s^k, s^k \rangle \qquad \forall \tau \in [0, \bar \tau_k].
    \end{equation*}
\end{lemma}
\begin{proof}
    Let \(k\) be fixed and suppose, by contradiction, that for any \( \varepsilon > 0 \) there exists \( \tau_\varepsilon \in [0, \varepsilon] \) such that the point \( x_\varepsilon := \bar x^k + \tau_\varepsilon^2 d^k + \tau_\varepsilon s^k \) satisfies
    \[
        \fbe(x_\varepsilon) > \fbe(x^k) - \sigma \Vert r^k \Vert^2 + \frac{\mu}{2} \tau_\varepsilon^2 \langle B_k s^k, s^k \rangle.  
    \]
    Taking the limit for \( \varepsilon \downto 0 \), the strict continuity of \( \fbe \) yields
    \[
        \fbe(\bar x^k) \geq \fbe(x^k) - \sigma \Vert r^k \Vert^2 > \fbe(x^k) - \gamma \frac{1 - \gamma L_f}{2} \Vert r^k \Vert^2.  
    \]
    Here, the last inequality is strict since \( r^k \neq 0 \).
    However, this contradicts \cref{prop:fbe-decrease-minimizers}, and hence the original assertion must be false. 
    The claim follows directly.
\end{proof}\noindent
Under mild assumptions \cref{alg:curvilinear-panoc} converges to critical points.
The proof is similar to that of \cite[Theorem 5.6]{themelis_forward-backward_2018}.
\begin{proposition}[Criticality of limit points] \label{prop:curvilinear-criticality}
    Suppose that \cref{assump:fbe-basic,assump:weak-convexity} hold.
    Then the following statements hold for the iterates generated by \cref{alg:curvilinear-panoc}:
    \begin{propenum}
        \item $r^k \to 0$ square-summably, and all limit points of $\seq{x^k}$ and $\seq{\bar x^k}$ are critical points; in fact we have that $\omega(\seq{x^k}) = \omega(\seq{\bar x^k}) \subseteq \fix T_\gamma$;
        \item $\seq{\fbe(x^k)}$ converges to a (finite) value $\varphi_\star$, and if $\seq{x^k}$ is bounded, then so does $\seq{\varphi(\bar x^k)}$.
    \end{propenum}
\end{proposition}
\begin{proof}
    Since \( \langle B_k s^k, s^k \rangle \leq 0 \), the linesearch condition \eqref{eq:curvilinear-linesearch} implies \[
        \fbe(x^{k+1}) \leq \fbe(x^k) - \sigma \Vert r^k \Vert^2,
    \]
    which is the same as in \cite{themelis_forward-backward_2018}.
    The reasoning from the proof of \cite[Theorem 5.6]{themelis_forward-backward_2018}, based on a standard telescoping argument, is then easily verified to hold for \cref{alg:curvilinear-panoc}.
\end{proof}
We now move on to the second-order convergence results. 
In particular, the following theorem establishes that $\lim_{k \to \infty} \langle B_k s^k, s^k \rangle = 0$, which is qualitatively similar to the classical result \cite[Theorem 3.1]{more_use_1979} for smooth problems.
A notable point of difference, however, is that \cref{alg:curvilinear-panoc} does not require a level-boundedness assumption.
Only boundedness of the directions $\seq{d^k}$, $\seq{s^k}$ is assumed, which can always be explicitly enforced, if needed, without breaking the global convergence (cf.\,\cref{prop:curvilinear-criticality}) or the local superlinear rate (cf.\,\cref{th:curvi-superlinear}).
\begin{theorem} \label{th:curvi-second-order}
    Suppose that \cref{assump:fbe-basic,assump:weak-convexity} hold, and consider the iterates generated by \cref{alg:curvilinear-panoc}.
    If $\seq{d^k}$, $\seq{s^k}$ are bounded and if \cref{assump:fbe-twice-diff} holds with respect to the limit points of $\seq{\bar x^k}$, then \[
        \lim_{k \to \infty} \langle B_k s^k, s^k \rangle = 0.
    \]
\end{theorem}
\begin{proof}
    The sequence $\seq{\fbe(x^k)}$ decreases monotonically and is bounded below.
    Hence, $\lim_{k \to \infty} \fbe(x^k) - \fbe(x^{k+1}) = 0$.
    We distinguish two cases.
    First, if there exists a $\tau \in (0, 1]$ such that $\tau_k \geq \tau$ for all $k \in \N$, then we obtain from \eqref{eq:curvilinear-linesearch} that \[
        \sigma \Vert r^k \Vert^2 - \frac{\mu}{2} \tau^2 \langle B_k s^k, s^k \rangle \leq \fbe(x^k) - \fbe(x^{k+1}).
    \]
    By nonnegativity of $\sigma \Vert r^k \Vert^2$ and $- \frac{\mu}{2} \tau^2 \langle B_k s^k, s^k \rangle$, the conclusion follows immediately.
    Otherwise, we assume without loss of generality that $\lim_{k \to \infty} \tau_k = 0$.
    Then, for $k$ sufficiently large, we define $\sigma_k := \tau_k / \beta$ and $y^{k+1} := \bar x^k + \sigma_k^2 d^k + \sigma_k s^k$, for which we have by \eqref{eq:curvilinear-linesearch} that 
    \begin{equation} \label{prf:curvilinear-second-order-sigma}
        \fbe(y^{k+1}) > \fbe(x^k) - \sigma \Vert r^k \Vert^2 + \frac{\mu}{2} \sigma_k^2 \langle B_k s^k, s^k \rangle.
    \end{equation}
    Under \cref{assump:fbe-twice-diff} the FBE is of class \(\C^2\) around the limit points of \(\seq{\bar x^k}\).
    Hence, for $k$ sufficiently large, $\fbe$ is of class $\C^2$.
    By a Taylor series argument, we therefore obtain that \[
        \fbe(y^{k+1}) = \fbe(\bar x^k) + \langle \nabla \fbe(\bar x^k), \sigma_k^2 d^k + \sigma_k s^k \rangle + \frac{\sigma_k^2}{2} \langle \nabla^2 \fbe(\bar x^k) (\sigma_k^2 d^k + s^k), (\sigma_k^2 d^k + s^k)\rangle + o(\sigma_k^2).
    \]
    The fact that \(d^k\) and \( s^k\) are descent directions, i.e., $\langle \nabla \fbe(\bar x^k), d^k \rangle \leq 0$ and $\langle \nabla \fbe(\bar x^k), s^k \rangle \leq 0$, combined with the inequality \(
    \fbe(\bar x^k) \leq \fbe(x^k) - \sigma \Vert r^k \Vert^2 \) yield \[
        \fbe(y^{k+1}) \leq \fbe(x^k) - \sigma \Vert r^k \Vert^2 + \frac{\sigma_k^2}{2} \langle \nabla^2 \fbe(\bar x^k) s^k, s^k\rangle + o(\sigma_k^2).
    \]
    Combined with \eqref{prf:curvilinear-second-order-sigma}, this implies
    \[
        \lim_{k \to \infty} - \frac{1 - \mu}{2} \langle B_k s^k, s^k\rangle \leq \lim_{k \to \infty} \frac{1}{2} \langle (\nabla^2 \fbe(\bar x^k) - B_k) s^k, s^k \rangle = 0.
    \]
    Our desired conclusion follows from the nonnegativity of $-\langle B_k s^k, s^k \rangle$.
\end{proof}\noindent
By construction, the negative curvature directions are such that \(\langle B_k s^k, s^k \rangle \leq 0\), with equality holding only when \(\lambda_{\min}(B_k) \geq 0\).
The result \( \lim_{k \to \infty} \langle B_k s^k, s^k \rangle = 0 \) therefore implies that \( \lambda_{\min}(B_\star) = \lambda_{\min}(\nabla^2 \fbe(x^\star)) \geq 0 \), for any \( x^\star \in \omega(\seq{\bar x^k}) \subseteq \fix T_\gamma \), or that \(\lim_{k \to \infty} \Vert s^k \Vert \to 0\).
Well-chosen directions \( \seq{s^k} \) are such that the latter can only hold if \(\lambda_{\min}(B_\star) \geq 0\).
This can be achieved, for example, by selecting \(s^k\) such that for all \(\lambda_{\min}(B_k) < 0\) either \(\Vert s^k \Vert\) is constant, or scales proportional to \(-\lambda_{\min}(B_k)\). 
Various conditions on the sequence \( \seq{s^k} \) have been proposed in the literature. Here we present a straightforward modification of \cite[Eq. 1.4]{lucidi_curvilinear_1998}, originally introduced for negative curvature directions generated by a Lanczos method.
\begin{corollary}[Convergence to second-order stationary points]
    Suppose that the assumptions of \cref{th:curvi-second-order} are satisfied, and that the directions \( \seq{s^k} \) are chosen such that
    \begin{equation*}
        \begin{cases}
            \lim_{k \to \infty} \langle B_k s^k, s^k \rangle = 0\\
            \lim_{k \to \infty} \Vert \bar r^k \Vert = 0
        \end{cases} \qquad \text{implies} \qquad \liminf_{k \to \infty} \lambda_{\min}(B_k) \geq 0.
    \end{equation*}
    Then, the limit points of \( \seq{\bar x^k} \) are second-order stationary points of \eqref{eq:problem}.
\end{corollary}
Finally, we show that \cref{alg:curvilinear-panoc} converges at a local superlinear rate to strong local minima when Newton directions \( d^k \) are (eventually) used.
\begin{theorem}[Superlinear convergence] \label{th:curvi-superlinear}
    Suppose that \cref{assump:fbe-basic,assump:weak-convexity} hold, and that the iterates \( \seq{\bar x^k} \) converge to a strong local minimizer \( x^\star \), with respect to which \cref{assump:fbe-twice-diff} holds.
    If the directions \( d^k \) are (eventually) chosen as Newton-directions, i.e., \( d^k = d_N^k := - B_k^{-1} \nabla \fbe(\bar x^k) \) for all \( k \) sufficiently large, then eventually \( \tau_k = 1 \) is always accepted and the sequences $\seq{x^k}$, $\seq{\bar x^k}$, and $\seq{r^k}$ converge with a superlinear rate.
\end{theorem}
\begin{proof}
    For \( \bar x^k \) sufficiently close to \( x^\star \), and hence for sufficiently large \( k \), we have by \cref{corr:hessian-gn-connection} that \[ 
        \hat \partial^2 \fbe(x) = \left\{ B(x) := \gamma^{-1} Q_\gamma(x) (\id - P_\gamma(x) Q_\gamma(x)) \mid P_\gamma(x) = J \prox_{\gamma g}(x - \gamma \nabla f(x)) \right\} 
    \] becomes a singleton, with \( B \) continuous. 
    By \cite[Theorem 4.11]{themelis_forward-backward_2018} it follows that \( B_\star := B(x^\star) = \nabla^2 \fbe(x^\star) \succ 0 \), and hence we have that \( B_k = B(\bar x^k) \succ 0 \) for \( k \) sufficiently large.
    Consequently, we have \( s^k = 0 \) and the linesearch condition \eqref{eq:curvilinear-linesearch} reduces to
    \begin{equation}\label{prf:reduced-linesearch-condition}
        \fbe(\bar x^k + \tau_k^2 d^k) \leq \fbe(x^k) - \sigma \Vert r^k \Vert^2
    \end{equation}
    for large enough \( k \), which we now show is eventually accepted with \( \tau_k = 1 \) for the Newton direction \( d^k = d_N^k \).
    We have from a second-order Taylor expansion of \( \fbe \) around \( x^\star \) that
    \begin{equation*}
        \begin{aligned}
            \fbe(\bar x^k) &= \fbe(x^\star) + \frac{1}{2} \langle B_\star (\bar x^k - x^\star), \bar x^k - x^\star \rangle + o(\Vert \bar x^k - x^\star \Vert^2)\\
            \fbe(\bar x^k + d_N^k) &= \fbe(x^\star) + \frac{1}{2} \langle B_\star (\bar x^k + d_N^k - x^\star), \bar x^k + d_N^k - x^\star \rangle + o(\Vert \bar x^k + d_N^k - x^\star \Vert^2)
        \end{aligned}
    \end{equation*}
    By similar arguments as in the final part of the proof of \cref{th:quadratic-convergence}, we know that \( \Vert \bar x^k + d_N^k - x^\star \Vert = o(\Vert \bar x^k - x^\star \Vert) \) and hence we obtain
    \begin{equation*}
        \begin{aligned}
            \fbe(\bar x^k + d_N^k) - \fbe(\bar x^k) &= -\frac{1}{2} \langle B_\star (\bar x^k - x^\star), \bar x^k - x^\star \rangle + o(\Vert \bar x^k - x^\star \Vert^2)\\
            &\leq -\frac{1}{2} \lambda_{\min}(B_\star) \Vert \bar x^k - x^\star \Vert^2 + o(\Vert \bar x^k - x^\star \Vert^2).
        \end{aligned}
    \end{equation*}
    Therefore, \( \fbe(\bar x^k + d_N^k) \leq \fbe(\bar x^k) \) holds for \( k \) large enough, which implies that eventually,
    \begin{equation*}
        \fbe(\bar x^k + d_N^k) \leq \fbe(\bar x^k) \leq \fbe(x^k) - \gamma \frac{1 - \gamma L_f}{2} \Vert r^k \Vert^2 \leq \fbe(x^k) - \sigma \Vert r^k \Vert^2.
    \end{equation*}
    Hence, eventually \cref{prf:reduced-linesearch-condition} is satisfied for \( \tau_k = 1 \) with \( d^k = d_N^k \), meaning that the Newton steps are always accepted for sufficiently large \( k \).
    Thus, we have \(
        \lim_{k \to \infty} \Vert x^{k+1} - x^\star \Vert / \Vert \bar x^k - x^\star \Vert = 0.
    \)
    The sequence \( \seq{x^k} \) converges to \( x^\star \) by \cref{prop:curvilinear-criticality}, and its superlinear convergence follows from
    \begin{equation*}
        \begin{aligned}
            \Vert \bar x^k - x^\star \Vert &\leq \gamma \Vert R_\gamma(x^k) \Vert + \Vert x^k - x^\star \Vert = (\gamma L_R + 1) \Vert x^k - x^\star \Vert.
        \end{aligned}
    \end{equation*}
    Here, the first inequality uses the triangle inequality and the second uses the \( L_R \)-Lipschitz continuity of \( R_\gamma \) in a neighborhood around \( x^\star \) \cite[Theorem 4.7 (i)]{themelis_forward-backward_2018}.
    Then, the bound \( \Vert r^k \Vert \leq L_R \Vert x^k - x^\star \Vert \) directly yields the superlinear convergence of \( \seq{r^k} \), and likewise \( \Vert \bar x^k - x^\star \Vert \leq \gamma \Vert r^k \Vert + \Vert x^k - x^\star \Vert \) proves the superlinear convergence of \( \seq{\bar x^k} \).
\end{proof}\noindent
A similar result can be obtained when the positive curvature directions satisfy a milder Dennis-Mor\'e condition.
\begin{theorem}[Superlinear convergence under Dennis-Mor\'e condition]
    Suppose that the conditions from \cref{th:curvi-superlinear} are satisfied, but that the directions \( d^k \) only satisfy the Dennis-Mor\'e condition
    \begin{equation*} %
        \lim_{k \to \infty} \frac{\Vert Q_k \bar r^k + \nabla^2 \fbe(x^\star) d^k \Vert}{\Vert d^k \Vert} = 0.
    \end{equation*}
    Then, eventually the stepsize $\tau_k = 1$ is always accepted and the sequences $\seq{x^k}$, $\seq{\bar x^k}$, and $\seq{r^k}$ converge with a superlinear rate.
\end{theorem}
\begin{proof}
    Tracing the arguments of the proof of \cref{th:curvi-superlinear}, the linesearch condition \eqref{eq:curvilinear-linesearch} eventually reduces to \eqref{prf:reduced-linesearch-condition}.
    Therefore, we can repeat the arguments from the proof of \cite[Theorem 5.10]{themelis_forward-backward_2018}, where the only modifications involve (i) using local Lipschitz continuity of \( \nabla \fbe \) around \( x^\star \) (since \( \fbe \) is locally of class \( \C^2 \)) instead of local Lipschitz continuity of \( R_\gamma \); and (ii) using nonsingularity of \( \nabla^2 \fbe(x^\star) \) instead of \( J R_\gamma(x^\star) \). 
\end{proof}

\begin{remark}[Directions based on the fixed-point residual] \label{remark:curvilinear-path}
    Linesearch methods like \zerofpr{} and \panoc{} use directions \( d^k \) which approximate a Newton step on the fixed-point residual \( \bar r^k \), rather than on \(Q_k \bar r^k\).
    In practice, Broyden-type directions like L-BFGS are used, and a local superlinear convergence rate is established under a Dennis-Mor\'e condition of the form
    \begin{equation*}
        \lim_{k \to \infty} \frac{\Vert \bar r^k + J R_\gamma(x^\star) d^k \Vert}{\Vert d^k \Vert} = 0.
    \end{equation*}
    It is straightforward to show that \cref{alg:curvilinear-panoc} also obtains a local superlinear convergence rate for these directions, as the proof reduces to that of \cite[Theorem 5.10]{themelis_forward-backward_2018}.
    Yet, special care is needed regarding the convergence to second-order stationary points.
    Since in general $\langle Q_k \bar r^k, d^k \rangle \leq 0$ can no longer be guaranteed, the steps below \eqref{prf:curvilinear-second-order-sigma} may not hold.
    Possible solutions include (i) enforcing \( d^k \) to be a descent direction, or (ii) using an unconventional curvilinear path \( x(\tau) = \bar x + \tau^3 d + \tau s \) such that the term involving \( d \) vanishes faster as \( \tau \downto 0 \). 
    However, in our experience directions based on \( \nabla \fbe \) perform significantly better in practice, especially when second-order information of \( f \) is used anyway for the computation of \( \seq{s^k} \).
\end{remark}

%% file: numerics/main.tex
This section validates the practical applicability of the proposed second-order proximal gradient methods.
First, \cref{sec:numerics-motivating-examples} reconsiders the motivating examples from \cref{sec:related-work} and visualizes the iterates generated by \cref{alg:fbtr,alg:curvilinear-panoc}, illustrating how nonsmooth strict saddle points are avoided in which first-order methods may get trapped.
Then, \cref{sec:numerics-sparse-pca,sec:numerics-phase-retrieval} compare the proposed algorithms against \panoc{} \cite{stella_simple_2017} on the tasks of sparse principal component analysis (PCA) and phase retrieval.

All methods are initialized randomly unless stated otherwise, and all algorithms use the same random initial points in comparisons.
First-order methods are terminated when \( \Vert R_\gamma(x^k) \Vert_{\infty} \leq 10^{-10} \), and for the proposed second-order methods we additionally require that \( \lambda_{\min}(B_k) \geq - 10^{-10} \).
All methods determine the proximal gradient stepsize \( \gamma > 0 \) adaptively by checking the relevant quadratic upper bound \cite{de_marchi_proximal_2022}.
\Cref{alg:fbtr} uses Steihaug's conjugate gradient method \cite{steihaug_conjugate_1983} for solving the trust-region problems with tolerance \( \varepsilon_{CG} = \min \{ 0.5 \Vert \nabla \fbe(x^k) \Vert_{\infty}, \Vert \nabla \fbe(x^k) \Vert_{\infty}^{\nicefrac{3}{2}} \} \).
Moreover, the radius update rule \eqref{eq:radius-rule} is implemented with \( c_1 = 0.35, c_2 = 1, c_3 = 1.5 \) and \( \mu_1 = 0.5, \mu_2 = 0.7 \).
\Cref{alg:curvilinear-panoc} uses parameters \( \mu = 0.1, \beta = \nicefrac{1}{\sqrt{2}}, \sigma = \beta \gamma \frac{1 - \gamma L_f}{2} \).
\Cref{alg:curvilinear-panoc} and \panoc{} use L-BFGS with a buffer size of \(5\), and the former uses a Lanczos method for computing the directions of negative curvature.

\subsection{Motivating examples revisited} \label{sec:numerics-motivating-examples}

\begin{figure}[h]
    \begin{subfigure}{0.5\textwidth}
        \centering
        \includegraphics[width=\textwidth]{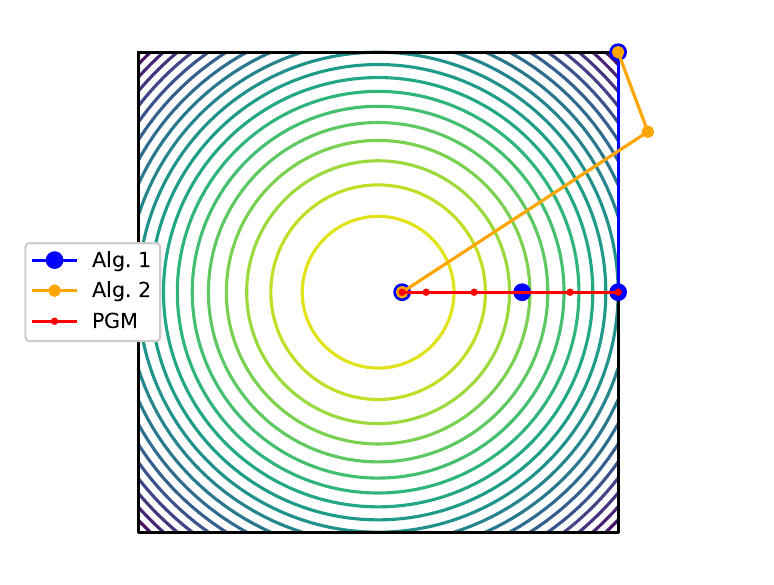}
        \captionsetup{width=1\textwidth}
        \caption[]{\centering$\varphi(x, y) = -x^2 - y^2 + \delta_{[-1, 1]}(x, y)$}
        \label{fig:quadratic-numerics}
    \end{subfigure}
    \begin{subfigure}{0.49\textwidth}
        \centering
        \includegraphics[width=\textwidth]{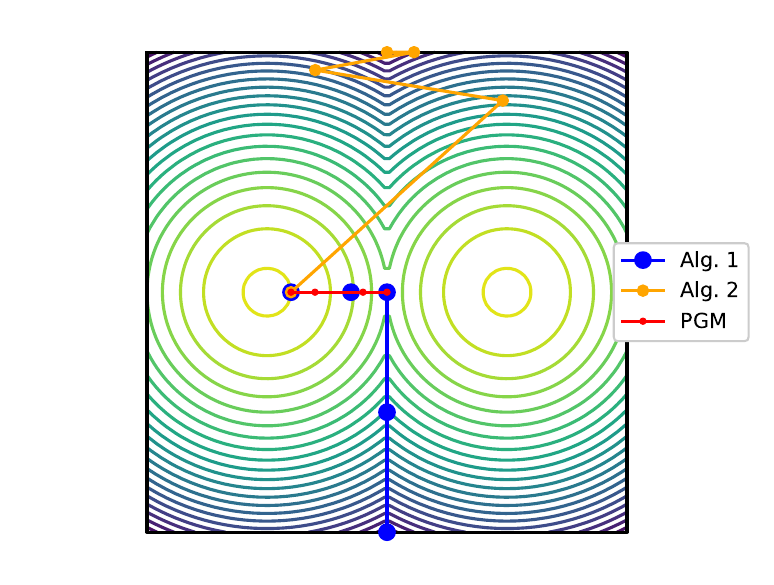}
        \captionsetup{width=.85\textwidth}
        \caption[]{$\varphi(x, y) = - x^2 - y^2 + \vert x \vert + \delta_{[-1, 1]}(x, y)$.}
        \label{fig:sparsity-numerics}
    \end{subfigure}
    \caption{Iterates of the PGM and \cref{alg:fbtr,alg:curvilinear-panoc} for the problems from \cref{ex:quadratics-box,ex:l1}, along with the level-curves of the corresponding objectives.}
    \label{fig:motivating-examples-numerics}
\end{figure}

Let us once more consider \cref{ex:quadratics-box}, which minimizes $\varphi(x, y) = -x^2 - y^2 + \delta_{[-1, 1]}(x, y)$.
Recall that the origin is a maximizer, that the points $(\pm 1, 0)$ and $(0, \pm 1)$ are nonsmooth strict saddle points, and that the points $(\pm 1, \pm 1)$ are strong local minimizers.
\Cref{fig:quadratic-numerics} visualizes the iterates generated by \cref{alg:fbtr,alg:curvilinear-panoc} and by the PGM, all initialized at $x^0 = (0.1, 0)$.
Observe how the PGM converges to the strict saddle point \( (1, 0) \), whereas both of the proposed methods converge to the strong local minimizer \( (1, 1) \).
Remark also how the iterates of \cref{alg:fbtr} pass through the strict saddle point but escape it, whereas \cref{alg:curvilinear-panoc} immediately moves towards \( (1, 1) \).

In a similar way, \cref{fig:sparsity-numerics} depicts the iterates of these algorithms when minimizing \( \varphi(x, y) = - x^2 - y^2 + \vert x \vert + \delta_{[-1, 1]}(x, y) \), which is the objective considered in \cref{ex:l1}, restricted to a compact box to ensure lower-boundedness.
We recall that the points \( (\pm 0.5, 0) \) are maximizers, that the points \( (0, 0) \) and \( (\pm 1, 0) \) are strict saddle points, and that the points \( (\pm 1, \pm 1) \) and \( (0, \pm 1) \) are strong local minimizers.
All methods are initialized at \( x^0 = (-0.4, 0) \).
As in the previous example, PGM converges to a strict saddle point.
In contrast, \cref{alg:fbtr,alg:curvilinear-panoc} converge to a strong local minimizer. 

\subsection{Sparse principal component analysis} \label{sec:numerics-sparse-pca}

Let us consider the sparse principal component analysis (PCA) problem described in \cite[\S 2.1]{journee_generalized_2010}, i.e.,
\begin{equation*}
    \minimize_{x \in \R^n} - \frac{1}{2} x^\top \Sigma x + \kappa \Vert x \Vert_1 + \delta_{\bar \cB(0; 1)}(x),
\end{equation*}
with sample covariance matrix \( \Sigma = A^\top A \), and sparsity inducing parameter \( \kappa \geq 0 \).
This can be formulated in the form \eqref{eq:problem} by defining the functions \( f(x) := - \frac{1}{2} x^\top \Sigma x \) and \( g(x) := \kappa \Vert x \Vert_1 + \delta_{\bar \cB(0; 1)}(x) \).
As in \cite[\S V]{themelis_new_2020} we generate sparse random matrices \( A \in \R^{20n \times n} \) with \( 10\% \) nonzeros.
The performance of \cref{alg:fbtr,alg:curvilinear-panoc} and \panoc{} is compared in \cref{table:sparse_pca} for \( n \in \{ 1000, 1500 \} \) and \( \kappa = 10^{-2} \) by counting the total number of matrix-vector products, as this constitutes the most costly operation involved.

Note that the proposed second-order methods significantly reduce the number of iterations, as well as the number of calls to the first-order oracles.
This confirms that, also performance-wise, the available second-order information is exploited. 
Remark that for \cref{alg:fbtr} the overall performance is roughly on par with that of \panoc{}.
As for \cref{alg:curvilinear-panoc}, we observe a considerably larger number of Hessian-vector products.
This is because (i) every iteration computes an (approximate) eigenvector \( s^k \) using a Lanczos procedure which on average performs more Hessian-vector products than Steihaug's CG method, and (ii) \cref{alg:curvilinear-panoc} performs significantly less iterations than \panoc{}, but still significantly more than \cref{alg:fbtr}.
This may be related to the question of effectively scaling the negative curvature direction, which trust-region methods do inherently.

\begin{table}[h]
    \caption{
        Comparison of \cref{alg:fbtr,alg:curvilinear-panoc}, and \panoc{} for sparse PCA problems.
        Median values over \( 100 \) random problem realizations are reported, with `hvp' denoting Hessian-vector products, and `mvp' matrix-vector products. 
    }
    \label{table:sparse_pca}
    \centering
    
    \begin{adjustbox}{width=\textwidth}
    \setlength\extrarowheight{3pt}
    \pgfplotstabletypeset[%
        begin table={\begin{tabular}[t]},
        every head row/.style={
            before row={%
              \hline
              \vphantom{\(\kappa\)}\\
              \hline
              \vphantom{ZeroFPR}\\
              \hline
            },
        },
        header=true,
        col sep=&,
        row sep=\\,
        string type,
        columns/{iters}/.style ={column name={iters}, column type={|l}},
        every row no 5/.style={after row=\hline},
        every row no 6/.style={after row=\hline},
    ]{
        \\
        iters\\
        eval \( f \)\\
        eval \( g \)\\
        grad \( f \)\\
        prox \( g \)\\
        Jprox \( g \)\\
        hvp \( f \)\\
        mvp\\
    }%
    \begin{filecontents*}{numerics/sparse_pca_1000_0.01_100.csv}
27.0	121.0	270.0
56.0	266.0	570.0
28.0	265.0	284.5
56.0	266.0	570.0
28.0	265.0	284.5
27.0	121.0	0.0
508.0	7620.0	0.0
564.0	7886.0	570.0
\end{filecontents*}
\pgfplotstabletypeset[%
        begin table={\begin{tabular}[t]},
        every head row/.style={
        before row={%
          \hline
          \multicolumn{3}{|c|}{\(n = 1000\)}\\
          \hline
        },
        after row/.add={}{\hline},
        },
        header=true,
       precision=1,
       columns/0/.style ={column name={\cref{alg:fbtr}}, column type={|l}},
       columns/1/.style ={column name={\cref{alg:curvilinear-panoc}}, column type={c}},
       columns/2/.style ={column name={PANOC}, column type={r|}},
        every row no 7/.style={after row=\hline},
        every row no 6/.style={after row=\hline},
    ]{numerics/sparse_pca_1000_0.01_100.csv}%
    \begin{filecontents*}{numerics/sparse_pca_1500_0.01_100.csv}
33.0	137.5	300.5
68.0	296.0	636.0
34.0	295.0	317.5
68.0	296.0	636.0
34.0	295.0	317.5
33.0	137.5	0.0
625.0	8659.5	0.0
693.0	8955.5	636.0
\end{filecontents*}
\pgfplotstabletypeset[%
        begin table={\begin{tabular}[t]},
        every head row/.style={
        before row={%
          \hline
          \multicolumn{3}{c|}{\(n = 1500\)}\\
          \hline
        },
        after row/.add={}{\hline},
        },
        header=true,
       precision=1,
       columns/0/.style ={column name={\cref{alg:fbtr}}, column type={l}},
       columns/1/.style ={column name={\cref{alg:curvilinear-panoc}}, column type={c}},
       columns/2/.style ={column name={PANOC}, column type={r|}},
        every row no 7/.style={after row=\hline},
        every row no 6/.style={after row=\hline},
    ]{numerics/sparse_pca_1500_0.01_100.csv}
    \end{adjustbox}
    
\end{table}

\subsection{Phase retrieval} \label{sec:numerics-phase-retrieval}

Let us consider the real-valued phase retrieval problem described in  \cite{sun_geometric_2018}, i.e.,
\begin{equation*}
    \minimize_{x \in \R^n} \frac{1}{2 m} \sum_{i = 1}^{m} (y_i^2 - (a_i^\top x)^2)^2 + \delta_{\bar \cB(0;1)}(x)
\end{equation*}
where \( a_i \in \R^n, i \in \N_{[1, N]} \) are the measurement vectors, and \( y_i \in \R, i \in \N_{[1, N]} \) the corresponding measurements.
The nonsmooth term \( \delta_{\bar \cB(0; 1)} \) is added to further exploit the known magnitude, which is assumed \( 1 \) without loss of generality.
Problems are constructed by generating standard Gaussian measurement vectors \( a_i \in \R^n, i \in \N_{[1, N]} \) and a random normalized solution \( x^\star \).
The noiseless measurements \( y_i \in \R, i \in \N_{[1, N]} \) are computed from these.
Henceforth, we fix \( n = 100 \) and consider \( m \in \{ 300, 3000 \} \).

In our experiments, we observed that for \( m = 3000 \) all methods consistently converged to the global optimal \( \varphi^\star = 0 \), whereas for \( m = 300 \) this was not the case.
This appears related to the benign landscape of the (unconstrained) phase retrieval problem when the number of measurements \( m \) is sufficiently large compared to the variable size \( n \) \cite{sun_geometric_2018}.
\Cref{table:phase_retrieval_best_global} compares the quality of the solutions obtained by \cref{alg:fbtr,alg:curvilinear-panoc} and \panoc{} for \(m = 300\).
It is remarkable that \cref{alg:curvilinear-panoc} finds significantly better solutions compared to both \cref{alg:fbtr} and \panoc{}.
Moreover, it turns out that the scaling of the negative curvature direction \( s^k \) in \cref{alg:curvilinear-panoc} not only has an impact on the overall performance, but also on the quality of the solution.
Finally, \cref{table:phase_retrieval} details the median number of calls to the first and second-order oracles by each of the methods, with similar conclusions as in the sparse PCA example.
In particular, \cref{alg:fbtr} is again roughly on par with \panoc{}, whereas \cref{alg:curvilinear-panoc} performs a large number of Hessian-vector products.

\begin{table}[h]
    \caption{
        Number of times each solver attained, up to a tolerance of \( 10^{-3} \), the best found objective and the global optimal \( \varphi^\star = 0 \), respectively, in solving \( 100 \) phase retrieval problems with \(n = 100, m = 300\).
    }
    \label{table:phase_retrieval_best_global}
    \centering
    
    \setlength\extrarowheight{3pt}
    \pgfplotstabletypeset[%
        begin table={\begin{tabular}[t]},
        every head row/.style={
            before row={%
            \hline
            \vphantom{\cref{alg:fbtr}}\\
            \vphantom{\( \bar s = 10^{-2} \)}\\
            \hline
            },
        },
        header=true,
        col sep=&,
        row sep=\\,
        string type,
        columns/{iters}/.style ={column name={Best objective found}, column type={|l}},
        every row no 0/.style={after row=\hline},
        every row no 6/.style={after row=\hline},
    ]{
        \\
        iters\\
        Global optimal found\\
    }%
    \begin{filecontents*}{numerics/phase_retrieval_best_global_100_300_100.csv}
11.0	96.0	29.0	18.0	14.0
10.0	80.0	27.0	16.0	12.0
\end{filecontents*}
\pgfplotstabletypeset[%
        begin table={\begin{tabular}[t]},
        every head row/.style={
        before row={%
          \hline
          \cref{alg:fbtr} & \multicolumn{3}{c|}{\cref{alg:curvilinear-panoc}} & PANOC \\
        },
        after row/.add={}{\hline},
        },
        header=true,
       precision=1,
       columns/0/.style ={column name={}, column type={|c}},
       columns/1/.style ={column name={\( \bar s = 1 \)}, column type={|c}},
       columns/2/.style ={column name={\( \bar s = 10^{-2} \)}, column type={c}},
       columns/3/.style ={column name={\( \bar s = 10^{-4} \)}, column type={c|}},
       columns/4/.style ={column name={}, column type={c|}},
       every row no 1/.style={after row=\hline},
    ]{numerics/phase_retrieval_best_global_100_300_100.csv}%
    
\end{table}

\vspace{-4mm}
\begin{table}[h]
    \caption{
        Comparison of \cref{alg:fbtr,alg:curvilinear-panoc} (with \( \bar s = 1 \)), and \panoc{} for solving phase retrieval problems with \(n = 100\).
        Median values over \( 100 \) random problem realizations are reported, with `hvp' denoting Hessian-vector products. 
    }
    \label{table:phase_retrieval}
    \centering
    
    \begin{adjustbox}{width=\textwidth}
    \setlength\extrarowheight{3pt}
    \pgfplotstabletypeset[%
        begin table={\begin{tabular}[t]},
        every head row/.style={
            before row={%
              \hline
              \vphantom{\(\kappa\)}\\
              \hline
              \vphantom{ZeroFPR}\\
              \hline
            },
        },
        header=true,
        col sep=&,
        row sep=\\,
        string type,
        columns/{iters}/.style ={column name={iters}, column type={|l}},
        every row no 5/.style={after row=\hline},
        every row no 6/.style={after row=\hline},
    ]{
        \\
        iters\\
        eval \( f \)\\
        eval \( g \)\\
        grad \( f \)\\
        prox \( g \)\\
        Jprox \( g \)\\
        hvp \( f \)\\
    }%
    \begin{filecontents*}{numerics/phase_retrieval_100_300_100.csv}
32.0	123.0	497.5
74.0	299.0	1250.0
37.0	298.0	624.5
74.0	299.0	1250.0
37.0	298.0	624.5
32.0	123.0	0.0
591.0	14899.5	0.0
\end{filecontents*}
\pgfplotstabletypeset[%
        begin table={\begin{tabular}[t]},
        every head row/.style={
        before row={%
          \hline
          \multicolumn{3}{|c|}{\(m = 300\)}\\
          \hline
        },
        after row/.add={}{\hline},
        },
        header=true,
       precision=1,
       columns/0/.style ={column name={\cref{alg:fbtr}}, column type={|l}},
       columns/1/.style ={column name={\cref{alg:curvilinear-panoc}}, column type={c}},
       columns/2/.style ={column name={PANOC}, column type={r|}},
        every row no 7/.style={after row=\hline},
        every row no 6/.style={after row=\hline},
    ]{numerics/phase_retrieval_100_300_100.csv}%
    \begin{filecontents*}{numerics/phase_retrieval_100_3000_100.csv}
14.0	29.0	54.5
32.0	69.0	145.0
16.0	68.0	72.0
32.0	69.0	145.0
16.0	68.0	72.0
14.0	29.0	0.0
130.0	1960.5	0.0
\end{filecontents*}
\pgfplotstabletypeset[%
        begin table={\begin{tabular}[t]},
        every head row/.style={
        before row={%
          \hline
          \multicolumn{3}{c|}{\(m = 3000\)}\\
          \hline
        },
        after row/.add={}{\hline},
        },
        header=true,
       precision=1,
       columns/0/.style ={column name={\cref{alg:fbtr}}, column type={l}},
       columns/1/.style ={column name={\cref{alg:curvilinear-panoc}}, column type={c}},
       columns/2/.style ={column name={PANOC}, column type={r|}},
        every row no 7/.style={after row=\hline},
        every row no 6/.style={after row=\hline},
    ]{numerics/phase_retrieval_100_3000_100.csv}%
    \end{adjustbox}
    
\end{table}

%% file: conclusion/main.tex
This work introduced second-order algorithms converging to second-order stationary points of a nonsmooth objective \( \varphi = f + g \), irrespective of the initial iterate.
To this end, we first established conditions under which the forward-backward envelope $\fbe$ is of class \( \C^2 \) locally around critical points, as well as an equivalence between the second-order stationary points of $\varphi$ and $\fbe$.
In particular, these results were shown to hold for the class of $\C^2$-partly smooth functions $g$ under a strict complementarity condition \eqref{eq:strict-complementarity}, thus demonstrating their relevance to a broad class of nonsmooth problems. 
We then presented and analyzed \cref{alg:fbtr}, a nonsmooth trust-region method, and \cref{alg:curvilinear-panoc}, a curvilinear linesearch method.
To the best of our knowledge, these two second-order methods are the first to provably escape strict saddle points of the composite problem \eqref{eq:problem}.

Numerical experiments validated these theoretical findings, and showcased the practical applicability of the proposed methods on the tasks of sparse principal component analysis and phase retrieval.
We observed that for these problems \cref{alg:fbtr} performs roughly on par with \panoc{}, an L-BFGS accelerated proximal gradient method.
In contrast, \cref{alg:curvilinear-panoc} performed subpar due to the negative curvature direction \( s^k \) being (i) expensive to compute, and (ii) hard to scale relative to the `fast' direction \( d^k \).
Yet, this method was found to converge to better stationary points, depending on the scaling of the negative curvature directions.

%% file: appendix/global-convergence.tex
This section proves, under conventional assumptions, the global convergence of \cref{alg:fbtr} when $\fbe \in \C^{1+}$, which is the case when \cref{assump:fbe-basic-tr} is satisfied (cf.\, \cref{prop:fbe-locally-c1}).
We also assume level boundedness of \(\fbe\), as guaranteed by \cref{lem:bounded-level-sets-fbe} under \cref{ass:lowerLevelSet}.
In particular, the sequence of gradients converges to zero, which corresponds to global subsequential convergence of the iterates $\seq{x^k}$.

The following condition is satisfied by \cref{lem:bounded-Hessian-approximants}.
\begin{condition}[Bounded Hessian approximants]\label{ass:lowerboundedB}
    The approximate Hessians $B_k$ are uniformly bounded, i.e.,
    \begin{equation*}
        \exists M > 0 : \forall k \in \N : \Vert B_k \Vert \leq M
    \end{equation*}
\end{condition}\noindent 
\begin{lemma}\label{lem:difference}
    Suppose that \cref{ass:lowerboundedB} holds. For the sequence $\seq{x^k}$ generated by algorithm \ref{alg:fbtr} we have
    \begin{equation*}\label{eq:nomRatio}
        \left|\fbe(x^k)-\fbe(x^{k} + d^k) -(m_k(0)-m_k(d^k))\right| = o(\| \delta_k \|)
    \end{equation*}
\end{lemma}
\begin{proof}  
    Continuous differentiability of the FBE means, by definition, that
    \begin{equation*}
        \lim_{d^k \to 0} \frac{\fbe(x^k + d^k) - \fbe(x^k) - \langle \nabla \fbe(x^k), d^k \rangle}{\Vert d^k \Vert} = 0.
    \end{equation*}
    Using little-$o$ notation, this is equivalent to
    $$
        \fbe (x^k + d^k) - \fbe(x^k) - \langle \nabla \fbe(x^k), d^k \rangle = o(\| d^k \|).
    $$
    Under Assumption \ref{ass:lowerboundedB} and by definition of $m_k$ we therefore have
    \begin{align*}
        \big|\fbe(x^k) -\fbe(x^k + d^k) &-(m_k(0)-m_k(d^k))\big|\\
        &= \left| m_k(d^k) - \fbe(x^k + d^k) \right|\\
        &= \left| \fbe(x^k) + \langle \nabla \fbe(x^k), d^k \rangle + \frac{1}{2} \langle B_k d^k, d^k \rangle - \fbe(x^k + d^k) \right|\\
        &\leq \left| \fbe (x^k + d^k) - \fbe(x^k) - \langle \nabla \fbe(x^k), d^k \rangle \right| + \frac{1}{2} \left| \langle B_k d^k, d^k\rangle \right| \\
        &\leq o(\|d^k\|)+\frac{1}{2}M \|d^k\|^2 \leq o(\delta_k).
    \end{align*}
\end{proof}\noindent
The following Theorem guarantees the \emph{well-definedness} of \cref{alg:fbtr}.
\begin{theorem} \label{lem:fbtr-succ-iterate}
    Suppose that \cref{ass:lowerboundedB,assump:suff-decrease} hold, that the sequence $\left(x^k\right)_{k \in \N}$ is generated by \cref{alg:fbtr} and that $\Vert \nabla \fbe(x^k) \Vert \geq \epsilon > 0$. Then, for any $k$, there exists a nonnegative integer $p \geq 0$ such that $x^{k+p+1}$ is a very successful iterate.
\end{theorem}

\begin{proof}
    Suppose that the statement does not hold for a given $k$, such that for any arbitrary $p \geq 0$, the iterate $x^{k+p+1}$ is not a very successful point. This means that for $p = 0, 1\dots$ we have that $\trratio_k < \mu_2$. The radius update rule \eqref{eq:radius-rule} therefore implies
    \begin{equation} \label{eq:fbtr-inf-succ-step1}
        \lim_{p \to \infty} \delta_{k+p} = 0.
    \end{equation}
    Moreover, we know from \cref{ass:lowerboundedB,assump:suff-decrease} and $\Vert \nabla \fbe(x^k) \Vert \geq \epsilon > 0$ that
    \begin{equation} \label{eq:fbtr-inf-succ-step2}
        m_k(0_n) - m_k(d^k) \geq \beta \Vert \nabla \fbe(x^k) \Vert \min \left\{ \delta_k, \frac{\Vert \nabla \fbe(x^k) \Vert}{\Vert B_k \Vert} \right\} \geq \beta \epsilon \min \left\{ \delta_k, \frac{\epsilon}{M} \right\}
    \end{equation}
    Thus, it follows from Lemma \ref{lem:difference} and equations \eqref{eq:fbtr-inf-succ-step1} and \eqref{eq:fbtr-inf-succ-step2} that
    \begin{equation*}
        \begin{aligned}
            \lvert \trratio_{k+p} - 1 \rvert &\overset{\eqref{eq:tr-ratio}}{=} \bigg\lvert \frac{\fbe(x^{k+p}) - \fbe(x^{k+p} + d^{k+p})}{m_{k+p}(0_n) - m_{k+p}(d^{k+p})} - 1\bigg\rvert \\
            &= \bigg\lvert \frac{\fbe(x^{k+p}) - \fbe(x^{k+p} + d^{k+p}) - (m_{k+p}(0_n) - m_{k+p}(d^{k+p}))}{m_{k+p}(0_n) - m_{k+p}(d^{k+p})}\bigg\rvert\\
            &\leq \frac{o(\Vert d^{k+p} \Vert)}{\beta \epsilon \min \left\{ \delta_{k+p}, \frac{\epsilon}{M} \right\}} \leq \frac{o(\delta_{k+p})}{\beta \epsilon \min \left\{ \delta_{k+p}, \frac{\epsilon}{M} \right\}} \to 0 ~ \text{as} ~ p \to \infty.
        \end{aligned}
    \end{equation*}
    Hence, for sufficiently large values of $p$, we have a very successful iterate $\trratio_{k+p} \geq \mu_2$. This contradiction completes the proof.
\end{proof}\noindent
The following theorem guarantees that there exists a subsequence of the gradients that converges to zero.
\begin{theorem} \label{th:fbe-limit-point}
    Suppose that \cref{assump:fbe-basic-tr,ass:lowerLevelSet} and \cref{ass:lowerboundedB,assump:suff-decrease} hold, then the sequence $\left( x^k \right)_{k \in \N}$ generated by \cref{alg:fbtr} satisfies
    \begin{equation*} \label{eq:fbe-limit-point}
        \liminf_{k \to \infty} \Vert \nabla \fbe (x^k) \Vert = 0.
    \end{equation*}
\end{theorem}

\begin{proof}
    This proof follows ideas from the proof of \cite[Theorem 4.5]{nocedal_numerical_2006}.
    \textit{Suppose for contradiction} that there exists an $\epsilon > 0$ and an index $K >0$ such that
    \begin{equation} \label{prf:fbe-limit-point-assumption}
        \Vert \nabla \fbe(x^k) \Vert \geq \epsilon \quad \quad \forall k \geq K.
    \end{equation}
    We know from \cref{ass:lowerboundedB,assump:suff-decrease} and $\Vert \nabla \fbe(x^k) \Vert \geq \epsilon > 0$ that for $k \geq K$:
    \begin{equation} \label{prf:fbe-limit-point-model-reduction}
        m_k(0_n) - m_k(d^k) \geq \beta \Vert \nabla \fbe(x^k) \Vert \min \left\{ \delta_k, \frac{\Vert \nabla \fbe(x^k) \Vert}{\Vert B_k \Vert} \right\} \geq \beta \epsilon \min \left\{ \delta_k, \frac{\epsilon}{M} \right\}.
    \end{equation}
    Then it follows from \cref{lem:difference}, \eqref{prf:fbe-limit-point-model-reduction} and $\Vert d^k \Vert \leq \delta_k$ that for $k \geq K$:
    \begin{equation} \label{prf:fbe-limit-point-ratio-minus-one}
        \begin{aligned}
            \lvert \trratio_{k} - 1 \rvert &\overset{\eqref{eq:tr-ratio}}{=} \bigg\lvert \frac{\fbe(x^{k}) - \fbe(x^{k} + d^{k})}{m_{k}(0_n) - m_{k}(d^{k})} - 1\bigg\rvert \\
            &= \bigg\lvert \frac{\fbe(x^{k}) - \fbe(x^{k} + d^{k}) - (m_{k}(0_n) - m_{k}(d^{k}))}{m_{k}(0_n) - m_{k}(d^{k})}\bigg\rvert\\
            &\leq \frac{o(\delta_k)}{\beta \epsilon \min \left\{ \delta_{k}, \frac{\epsilon}{M} \right\}}.
        \end{aligned}
    \end{equation}
    As \(\delta_k \to 0\), also \(\rho_k \to 1\).
    Consequently, for sufficiently small \(\delta_k > 0\), eventually \(\rho_k \geq \mu_2\) holds which by the update rule \eqref{eq:radius-rule} implies \(\delta_{k+1} \geq \delta_k\).
    We conclude that \(\delta_k\) is lower bounded.
    By \cref{lem:fbtr-succ-iterate}, there exists an infinite subsequence $\mathcal{K}$ such that $\trratio_k \geq \mu_1, \forall k \in \mathcal{K}$.
    Then it follows from \eqref{prf:fbe-limit-point-model-reduction} that for $k \in \mathcal{K}$ and $k \geq K$:
    \begin{equation*}
        \begin{aligned}
            \fbe(x^k) - \fbe(x^{k+1}) &= \fbe(x^k) - \fbe(x^k + d^k)\\
            &\overset{\trratio_k \geq \mu_1}{\geq} \mu_1 (m_k(0) - m_k(d^k))\\
            &\geq \mu_1 \beta \Vert \nabla \fbe(x^k) \Vert \min \left\{ \delta_k, \frac{\Vert \nabla \fbe (x^k)\Vert}{\Vert B_k \Vert} \right\} \\ 
            &\geq \mu_1 \beta \epsilon \min \left\{ \delta_k, \frac{\epsilon}{M} \right\}
        \end{aligned}
    \end{equation*}
    Since $\fbe$ is bounded below, this implies that
    \(
        \lim_{\substack{k \to \infty\\ k \in \mathcal{K}}} \delta_k = 0.
    \)
    However, this contradicts the lower boundedness of \(\delta_k\), 
    meaning that the original assertion \eqref{prf:fbe-limit-point-assumption} must be false.
    This completes the proof.
\end{proof}\noindent
Finally, we present a proof for \cref{th:fbe-limit}, showing that not only a subsequence, but the full sequence of gradients converges to zero.
\begin{theorem}[Global convergence]
    Suppose that \cref{assump:fbe-basic-tr,ass:lowerLevelSet} and \cref{ass:lowerboundedB,assump:suff-decrease} hold, then the sequence $\left( x^k \right)_{k \in \N}$ generated by \cref{alg:fbtr} satisfies
    \begin{equation*}
        \lim_{k \to \infty} \Vert \nabla \fbe (x^k) \Vert = 0
    \end{equation*}
\end{theorem}

\begin{proof}
    This proof closely follows the proof of \cite[Theorem 4.6]{nocedal_numerical_2006}, which is in turn based on that of \cite[Theorem 2.2]{shultz_family_1985}.
    Under Assumption \ref{ass:lowerLevelSet} the lower level-set $\mathcal{L}_\gamma(x^0)$ is bounded.
    By \cref{prop:fbe-C1} $\nabla \fbe$ is L-Lipschitz continuous on the bounded set $\mathcal{L}_\gamma(x^0)$, for some $L > 0$.
    Consider an arbitrary index $m$ such that $\nabla \fbe(x^m) \neq 0$.
    Then we have for all $x \in \mathcal{L}_\gamma(x^0)$ that
    \begin{equation*}
        \Vert \nabla \fbe(x) - \nabla \fbe(x^m) \Vert \leq L \Vert x - x^m \Vert
    \end{equation*}
    Define the scalars
    \begin{equation*}
        \epsilon := \frac{1}{2} \Vert \nabla \fbe(x^m) \Vert, \quad \quad R :=  \frac{\epsilon}{L}
    \end{equation*}
    and observe that the ball
    \begin{equation*}
        \mathcal{B}(x^m, R) = \left\{ x \mid \Vert x - x^m \Vert \leq R \right\}
    \end{equation*}
    is contained in $\mathcal{L}_\gamma(x^0)$ such that $\fbe$ is \(L\)-smooth in $\mathcal{B}(x^m, R)$.
    Moreover, note that
    \begin{equation*}
        x \in \mathcal{B}(x^m, R) \Rightarrow \Vert \nabla \fbe(x) \Vert \geq \Vert \nabla \fbe(x^m) \Vert - \Vert \nabla \fbe(x) - \nabla \fbe(x^m) \Vert \geq \epsilon.
    \end{equation*}
    Thus, if the entire sequence $\left\{ x^k \right\}_{k \geq m}$ stays inside the ball $\mathcal{B}(x^m, R)$, we would have that $\Vert \nabla \fbe(x^k) \Vert \geq \epsilon > 0$, for all $k \geq m$.
    By the same reasoning as in the proof of Theorem \ref{th:fbe-limit-point}, it can be shown that this scenario does not occur.
    Consequently, the sequence $\left\{ x^k \right\}_{k \geq m}$ eventually leaves the ball $\mathcal{B}(x^m, R)$.
    Define the index $\ell \geq m$ such that $x^{\ell + 1}$ is the first iterate after $x^{m}$ outside $\mathcal{B}(x^m, R)$.
    We can use $\Vert \nabla \fbe(x^k) \Vert \geq \epsilon$ for $k = m, \dots, \ell$ and \eqref{prf:fbe-limit-point-model-reduction} to obtain
    \begin{equation*}
        \begin{aligned}
            \fbe(x^m) - \fbe(x^{\ell+1}) &= \sum_{k = m}^\ell \fbe(x^k) - \fbe(x^{k+1})\\
            &\geq \sum_{\substack{k = m\\ x^k \neq x^{k+1}}}^\ell \mu_1 \left[ m_k(0_n) - m_k(d^k) \right]
            \geq \sum_{\substack{k = m\\ x^k \neq x^{k+1}}}^\ell \mu_1 \beta \epsilon \min \left\{ \delta_k, \frac{\epsilon}{M} \right\},
        \end{aligned} 
    \end{equation*}
    in which we limited the summation to those iterations where the candidate step was accepted, i.e., $\trratio_k \geq \mu_1$.
    We distinguish two cases.
    If $\delta_k \leq \nicefrac{\epsilon}{M}$ for all $k = m, \dots \ell$, then we have that
    \begin{equation*}
        \fbe(x^m) - \fbe(x^{\ell+1}) \geq \mu_1 \beta \epsilon \sum_{\substack{k = m\\ x^k \neq x^{k+1}}}^\ell \delta_k \geq \mu_1 \beta \epsilon R = \mu_1 \beta \epsilon \frac{\epsilon}{L}, 
    \end{equation*}
    where we made use of the fact that
    \begin{equation*}
        \begin{aligned}
            \sum_{\substack{k = m\\ x^k \neq x^{k+1}}}^\ell \delta_k &= \sum_{\substack{k = m\\ x^k \neq x^{k+1}}}^\ell \Vert x^{k+1} - x^k \Vert \geq \Vert x^{\ell+1} - x^{m} \Vert \geq R.
        \end{aligned}
    \end{equation*}
    Alternatively, $\delta_k > \nicefrac{\epsilon}{M}$ for at least some $k = m, \dots, \ell$, which implies
    \begin{equation*}
        \fbe(x^m) - \fbe(x^{\ell+1}) \geq \mu_1 \beta \epsilon \frac{\epsilon}{M}.
    \end{equation*}
   A combination of the above two cases leads to
    \begin{equation*}
        \begin{aligned}
            \fbe(x^m) - \fbe(x^{\ell+1}) &\geq \mu_1 \beta \epsilon \min \left\{ \frac{\epsilon}{L}, \frac{\epsilon}{M} \right\}\\
            &= \frac{1}{2} \mu_1 \beta \Vert \nabla \fbe(x^m) \Vert \min \left\{ \frac{\Vert \nabla \fbe(x^m) \Vert}{L}, \frac{\Vert \nabla \fbe(x^m) \Vert}{M} \right\} > 0.
        \end{aligned}
    \end{equation*}
    Since $\fbe$ is bounded below and $\{ \fbe(x^m) \}_{m \in \N}$ is monotonically decreasing, it follows that this sequence converges to some $\fbe^\star$ as $m \to \infty$.
    Therefore, as $m \to \infty$ we have that
    \begin{equation*}
        \begin{aligned}
            0 \leftarrow \fbe(x^m) - \fbe^\star & \geq \fbe(x^m) - \fbe(x^{\ell(m)+1})\\
            & \geq \frac{1}{2} \mu_1 \beta \Vert \nabla \fbe(x^m) \Vert \min \left\{ \frac{\Vert \nabla \fbe(x^m) \Vert}{L}, \frac{\Vert \nabla \fbe(x^m) \Vert}{M} \right\},
        \end{aligned}
    \end{equation*}
    which means that $\Vert \nabla \fbe(x^m) \Vert \to 0$ must hold.
    This completes the proof.
\end{proof}